\documentclass[a4paper,11pt,reqno]{amsart}

\usepackage{hyperref}
\usepackage{amssymb}
\usepackage{graphics}
\usepackage{graphicx}
\usepackage{epsfig}
\usepackage{multicol}
\usepackage{xypic}
\usepackage{dsfont}

\usepackage{amsmath}

\newcommand{\note}[1]{\vspace{5 mm}\par \noindent
  \marginpar{\textsc{Note}}
  \framebox{\begin{minipage}[c]{0.95 \textwidth}
      \tt #1 \end{minipage}}\vspace{5 mm}\par}

\usepackage{amsmath,amsfonts,amssymb,amsthm}

\usepackage{fouridx}

\usepackage{eufrak}
\usepackage{graphics}
\usepackage{graphicx}
\usepackage{epsfig}
\usepackage{multicol}
\usepackage{xypic}
\usepackage{amsmath}

\newtheorem{proposition}{Proposition}[section]
\newtheorem{definition}{Definition}[section]
\newtheorem{lemma}[definition]{Lemma}
\newtheorem{remark}[definition]{Remark}

\newtheorem{algorithm}[definition]{Algorithm}

\newcommand{\lp}{\left(}
\newcommand{\rp}{\right)}
\newcommand{\lc}{\left\{}
\newcommand{\rc}{\right\}}
\newcommand{\der}{\partial}

\newcommand{\bra}{\langle}
\newcommand{\ket}{\rangle}
\newcommand{\R}{\mathds{R}}      
\newcommand{\F}{\mathds{F}}
\newcommand{\I}{\mathds{I}}
\newcommand{\Flder}{\rightarrow}

\usepackage{amssymb}

\newcommand{\proa}{A^*G \mbox{$\;$}_{\tau^*} \kern-3pt\times_\alpha
G \mbox{$\;$}_\beta \kern-3pt\times_{\tau^*} A^*G}

\newcommand{\alg}{\mathfrak{so}(3)}

\newcommand{\e}{\mbox{exp}}
\newcommand{\Ad}{\mbox{Ad}}
\newcommand{\ca}{\mbox{cay}}
\newcommand{\ad}{\mbox{ad}}

\newcommand{\dmuno}{(\mbox{d}\tau^{-1}_{h\,\xi_{k}})^{*}}

\newcommand{\al}{\mathfrak{g}}
\newcommand{\dal}{\mathfrak{g}^{*}}

\newcommand{\setr}{\mathfrak{se}(3)}

\begin{document}

\title{Discrete Variational Optimal Control}

\author[F. Jim\'enez]{Fernando Jim\'enez}
\address{F. Jim\'enez: Instituto de Ciencias Matem\'aticas, CSIC-UAM-UC3M-UCM,
Campus de Cantoblanco, UAM,
C/Nicol\'as Cabrera, 15
 28049 Madrid, Spain} \email{fernando.jimenez@icmat.es}

\author[M. Kobilarov]{Marin Kobilarov}
\address{M. Kobilarov: California Institute of Technology, Control
  and Dynamical Systems, Pasadena, CA 91125, USA} \email{marin@cds.caltech.edu}

\author[D.\ Mart\'{\i}n de Diego]{David Mart\'{\i}n de Diego}
\address{D.\ Mart\'{\i}n de Diego: Instituto de Ciencias Matem\'aticas, CSIC-UAM-UC3M-UCM,
Campus de Cantoblanco, UAM,
C/Nicol\'as Cabrera, 15
 28049 Madrid, Spain} \email{david.martin@icmat.es}

\thanks{This work has been partially supported by MEC (Spain)
Grants   MTM 2010-21186-C02-01, MTM2009-08166-E,   and IRSES-project ``Geomech-246981''.}

\maketitle

\begin{abstract}

This paper develops numerical methods for optimal control of
mechanical systems in the Lagrangian setting. It extends the theory of
discrete mechanics to enable the solutions
of optimal control problems through the discretization of variational
principles. The key point is to solve the optimal control problem as a
variational integrator of a specially constructed higher-dimensional
system. The developed framework applies to systems on tangent bundles,
Lie groups, underactuated and nonholonomic systems with symmetries,
and can approximate either smooth or discontinuous control inputs. The
resulting methods inherit the preservation properties of variational
integrators and result in numerically robust and easily implementable
algorithms. Several theoretical and a practical examples, e.g. the control
of an underwater vehicle, will illustrate the application of the
proposed approach.

\end{abstract}

\section{Introduction}

The goal of this paper is to develop, from a geometric point of view,
numerical methods for optimal control of Lagrangian mechanical
systems. Our approach employs the theory of discrete mechanics and
variational integrators~\cite{mawest} to derive both an
integrator for the dynamics and an optimal control algorithm in a
unified manner. This is accomplished through the discretization of the
Lagrange-d'Alembert variational principle on manifolds. An integrator
for the mechanics is derived using a standard Lagrangian function and
virtual work done by control forces, while control optimality conditions
are derived using a special Lagrangian defined on a
higher-dimensional space which encodes the dynamics and a desired cost
function. The resulting integration and optimization schemes are
symplectic, respect the state space structure, and momentum
preserving. These qualities are associated with numerical stability
which motivate the development of practical algorithms that can be
applied to robotic or aerospace vehicles.

The proposed framework is general and applies to unconstrained systems, as
well systems with symmetries, underactuation, and nonholonomic
constraints. In particular, our construction is appropriate for
controlled Lagrangian systems that evolve on a general
tangent bundle $TQ$ with associated discrete state space $Q\times Q$,
where $Q$ is a differentiable manifold (\cite{mawest,objuma}). In
addition we focus on underactuated systems evolving on a Lie group $G$
(\cite{BlHuLeSa2009,BoSu,KM,LeMcLe2006}) that are applicable for
systems consisting of rigid bodies. Finally, the theory extends to
the more general principle bundle setting with discrete
analog $Q\times Q\times G$ (or more generally  $(Q\times Q)/G$)
assuming that the action of a Lie group $G$ of symmetry leaves the
control system invariant (\cite{Co2002,FeIgDe2008,KoMaSu2010}).

The main idea is the following: we take an approximation of the
Lagrange-d'Alembert principle for forced Lagrangian systems, which
models control inputs and external forces such as gravity or drag.
In principle, we admit the possibility of piecewise continuous control
forces, as happens in real applications.
We observe that the  discrete equations of motion for this type of
systems are interpreted as the discrete Euler-Lagrange equations of a
new Lagrangian defined in an augmented discrete phase space. Next, we
apply discrete variational calculus techniques to derive the discrete
optimality conditions. After this, we recover two sequences of
discrete controls  modeling a piecewise control trajectory.

Additionally we show how to derive the equations for different reduced systems.
We specifically develop numerical methods for systems on Lie groups
that lead to practical algorithm implementation. One such example
system--an underactuated underwater vehicle--is used to illustrate the
developed methodology. The resulting algorithm is simple to implement and
has the ability to quickly converge to a solution which is close to
the optimal solution and to the true system dynamics.
We also extend our techniques to more general reduced systems like
optimal control problems in trivial principal bundles and we show how
to introduce nonholonomic constraints in our framework.

\vspace{0.2cm}
Moreover, since we are reducing the optimality conditions to discrete Euler-Lagrange equations, the geometric preservation properties like sym\-plec\-tic-momentum preservation in the standard case or Poisson bracket and momentum preservation for reduced systems are automatically guaranteed  using the results in~\cite{groupoid,mawest}.
\vspace{0.2cm}

The paper is structured as follows: \S\ref{section2} introduces
variational integrators. \S\ref{TQ} formulates optimal control
problems for Lagrangian systems defined on tangent bundles, in the
continuous and discrete setting, and for both fully and underactuated
systems. A simple control problem for a mechanical Lagrangian on
$\R^n$ illustrates these developments. In \S\ref{RED}, discrete
mechanics on Lie groups is introduced. Specifically, discrete
Euler-Poincar\'e equations and their Hamiltonian version, the discrete
Lie-Poisson equations, are obtained. Sections \S\ref{DISC} and
\S\ref{Vehiculo} develop the discretization procedure and the
numerical aspects of the proposed approach. The developed algorithm is
illustrated with an application to an unmanned underwater vehicle
evolving on $SE(3)$. Finally, \S\ref{Ati} deals with reduced systems
on a trivial principal bundle and with nonholonomic mechanics.

\section{Discrete Mechanics and Variational Integrators}
\label{section2}
Let $Q$ be a $n$-dimensional differentiable manifold with local
coordinates $(q^i)$, $1\leq i\leq n$. Denote by $TQ$ its tangent
bundle with induced coordinates $(q^i, \dot{q}^i)$. Given a
Lagrangian function $L\colon TQ\to \R$ the Euler--Lagrange equations are
\begin{equation}\label{qwer}
\frac{d}{dt}\left(\frac{\partial L}{\partial\dot
q^i}\right)-\frac{\partial L}{\partial q^i}=0, \quad 1\leq i\leq
n.
\end{equation}
These equations are a system of implicit second order differential
equations. In the sequel, we will assume that the Lagrangian is
\textbf{regular}, that is, the matrix $\left(\frac{\partial^2
    L}{\partial \dot q^i \partial \dot q^j}\right)$ is non-singular.
It is well known that the origin of these equations is variational (see \cite{AbMa,MaRa}).

 Variational integrators retain  this variational character
and also some of main  geometric properties of the continuous
system, such as symplecticity and momentum conservation (see
\cite{Hair} and references therein).

In the following we will summarize the main features of this type
of numerical integrators \cite{mawest}.  A \textbf{discrete Lagrangian} is a map
$L_d:Q\times Q\rightarrow\R$, which may be considered as
an approximation of the integral action defined by a continuous  Lagrangian $L\colon TQ\to
\R$:
$
L_d(q_0, q_1)\approx \int^h_0 L(q(t), \dot{q}(t))\; dt
$
where $q(t)$ is a solution of the Euler-Lagrange equations for $L$, where $q(0)=q_0$ and $q(h)=q_1$ and $h>0$ is enough small.

\begin{remark}
{\rm
The Cartesian product $Q\times Q$ is equipped with an interesting differential structure, the Lie groupoid structure which allows us to extend the construction of variational calculus to another interesting situations (Lie groupoids). See \cite{groupoid} for more details.
}
\end{remark}

 Define the \textbf{action sum} $S_d\colon Q^{N+1}\to
\R$,   corresponding to the Lagrangian $L_d$ by
$
{S_d}=\sum_{k=1}^{N}  L_d(q_{k-1}, q_{k}),
$
where $q_k\in Q$ for $0\leq k\leq N$, and $N$ is the number of steps. The discrete variational
principle   states that the solutions of the discrete system
determined by $L_d$ must extremize the action sum given fixed
endpoints $q_0$ and $q_N$. By extremizing ${S_d}$ over $q_k$,
$1\leq k\leq N-1$, we obtain the system of difference equations
\begin{equation}\label{discreteeq}
 D_1L_d( q_k, q_{k+1})+D_2L_d( q_{k-1}, q_{k})=0,
\end{equation}
or, in coordinates,
\[
\frac{\partial L_d}{\partial x^i}(q_k, q_{k+1})+\frac{\partial
L_d}{\partial y^i}(q_{k-1}, q_{k})=0,
\]
where $1\leq i\leq n,\ 1\leq
k\leq N-1$ and $x,y$ denote the $n$-first and $n$-second variables of the function $L$ respectively.

These  equations are usually called the  \textbf{discrete
Euler--Lagrange equations}. Under some regularity hypotheses (the
matrix $(D_{12}L_d(q_k, q_{k+1}))$ is regular), it is possible to
define a (local) discrete flow $ \Upsilon_{L_d}\colon Q\times
Q\to  Q\times Q$, by $\Upsilon_{L_d}(q_{k-1}, q_k)=(q_k,
q_{k+1})$ from (\ref{discreteeq}). Define the  discrete
Legendre transformations associated to  $L_d$ as
\begin{eqnarray*}
\F^-L_d\colon Q\times Q&\to & T^*Q\\
(q_0, q_1)&\longmapsto& (q_0, -D_1 L_d(q_0, q_1)),\\
\F^+L_d\colon  Q\times Q&\to&  T^*Q\\
(q_0, q_1)&\longmapsto& (q_1, D_2 L_d(q_0, q_1))\; ,
\end{eqnarray*}
and the discrete Poincar{\'e}--Cartan 2-form $\omega_d=(\F^+L_d)^*\omega_Q=(\F^{-}L_d)^*\omega_Q$,
where $\omega_Q$ is the canonical symplectic form on $T^*Q$. The
discrete algorithm determined by $\Upsilon_{L_d}$ preserves the
symplectic form $\omega_d$, i.e., $\Upsilon_{L_d}^*\omega_d=\omega_d$.
Moreover, if the discrete Lagrangian is invariant under the
diagonal action of a Lie group $G$, then the discrete momentum map
$J_d\colon Q\times Q \to  {\mathfrak g}^*$ defined by
\[ \langle
J_d(q_k, q_{k+1}), \xi\rangle=\langle D_2L_d(q_k, q_{k+1}),
\xi_Q(q_{k+1})\rangle \]
is preserved by the discrete flow.
Therefore, these integrators are symplectic-momentum preserving. Here, $\xi_Q$ denotes the fundamental vector field
determined by $\xi\in {\mathfrak g}$, where ${\mathfrak g}$ is the Lie
algebra of $G$. (See \cite{mawest} for more details.)

\section{Discrete optimal control on tangent bundles}
\label{TQ}

Consider a mechanical system which configuration space is a
$n$-dimensional differentiable manifold $Q$ and which dynamics is
determined  by a Lagrangian $L: TQ\rightarrow \R$. The control forces
are modeled as a mapping
$f: TQ\times U\to T^*Q$, where $f(v_q, u)\in T_q^*Q$, $v_q\in T_qQ$
and $u\in U$, being $U$ the control space.
Observe that this last definition also covers configuration and
velocity dependent forces such as dissipation or friction (see \cite{objuma}).
For greater generality we consider control variables that are only
piecewise continuous to account for impulsive controls.

The motion of the mechanical system is described by applying the principle of {\bf Lagrange-D'Alembert}, which requires that the solutions $q(t)\in Q$ must satisfy
\begin{equation}\label{ldp1}
\delta \int^T_0 L(q(t), \dot{q}(t))\, dt + \int^T_0 f(q(t), \dot{q}(t), u(t))\,\delta q(t)\; dt=0,
\end{equation}
where $(q\,,\,\dot q)$ are the local coordinates of $TQ$ and where we consider arbitrary variations $\delta q\in T_{q(t)}Q$ with $\delta q(0)=0$ and $\delta q(T)=0$
(since we are prescribing fixed initial and final conditions $(q(0), \dot{q}(0))$ and $(q(T), \dot{q}(T))$).

Given that we are considering an optimal control problem, the forces $f$ must be chosen, if they exist, as the ones that extremize the \textbf{cost functional}:
\begin{equation}\label{cost-1}
\int_0^T C(q(t), \dot{q}(t), u(t))\, dt,
\end{equation}
where $C: TQ\times U\rightarrow \R$.

The optimal equations of motion can now be derived using Pontryagin
maximum principle. Generally, it is not possible to
explicitly integrate these equations and, consequently, it is
necessary to apply a numerical method. In this work, using discrete
variational techniques, we will first discretize the
Lagrange-d'Alembert principle and then the cost functional. We obtain
a numerical method that preserves some geometric features of the
original continuous system as we will see in the sequel.

To discretize this problem we replace the tangent space $TQ$ by the Cartesian product $Q\times Q$ and the continuous curves by sequences $q_0, q_1, \ldots q_N$ (we are using $N$ steps, with time step $h$ fixed, in such a way $t_{k}=kh$ and $Nh=T$).
The discrete Lagrangian $L_d: Q\times Q\rightarrow \R$ is constructed as an approximation of the action integral in a single time step (see \cite{mawest}), that is
\[
L_d(q_k, q_{k+1})\approx \int_{kh}^{(k+1)h}L(q(t), \dot{q}(t))\; dt.
\]
We choose the following discretization for the external forces: $f^{\pm}_{k}:Q\times Q\times U\Flder T^{*}Q$, where $U\subset\R^{m},\,\,m\leq n$, such that

\begin{eqnarray*}
f^-_{k}(q_k, q_{k+1}, u_{k}^{-})&\in&T^*_{q_k}Q,\\
 f^+_{k}(q_k, q_{k+1}, u_{k}^{+})&\in&T^*_{q_{k+1}}Q.
\end{eqnarray*}

Observe that, as mentioned above, we have introduced the discrete controls as two different sequences $\lc u_k^{-}\rc$ and $\lc u_k^{+}\rc$. In the notation followed through this paper, the time interval between $[k,k+1]$ is denoted as the $k$-th interval, while the controls in $k^+$ and $(k+1)^{-}$ are denoted by $u_{k}^{-}$ and $u_{k+1}^+$ respectively. This choice allows us to model piecewise continuous controls, admitting discrete jumps at the time steps $t_k=hk$. Our notation is completely depicted in the following figure:
\vspace{-0.5cm}

\begin{center}
\setlength{\unitlength}{1.8mm}
\begin{picture}(60,50)
  \put(0 ,15){\vector(1,0){59}}
  \put(60,14){$t_k$}
  \multiput(6,9)(15,0){4}{\line(0,1){35}}
  \put(5,5){$_{hk}$}
  \put(19,5){$_{h(k+1)}$}
  \put(34,5){$_{h(k+2)}$}
  \put(49,5){$_{h(k+3)}$}


  \put(31,17){$_{u_{k+1}^+}$}
  \put(36,17){\circle*{0.6}}
  \put(36,32){\circle*{0.6}}
  \put(36,32){\line(1,0){15}}
  \put(37,35){$_{u_{k+2}^{-}}$}

  \put(17,41){$_{u_k^{+}}$}
  \put(21,40){\circle*{0.6}}

  \put(22,32){$_{u_{k+1}^-}$}
  \put(21,32){\circle*{0.6}}
  \put(21,32){\line(1,-1){15}}

  \put(6,13){$\underbrace{\rule{2.7cm}{0cm}}$}
  \put(11,10){$_{(k)-th}$}

  \put(21,13){$\underbrace{\rule{2.7cm}{0cm}}$}
  \put(25,10){$_{(k+1)-th}$}

  \put(36,13){$\underbrace{\rule{2.7cm}{0cm}}$}
  \put(40,10){$_{(k+2)-th}$}

  \put(7,24){$_{u_k^{-}}$}
  \put(6,25){\circle*{0.6}}
  \put(6,25){\line(1,1){15}}

  \put(51,32){\circle*{0.6}}
  \put(46,35){$_{u_{k+2}^{+}}$}
\end{picture}
\end{center}
Moreover, we have that
\begin{eqnarray*}
  &&f^-_{k}(q_k, q_{k+1}, u_{k}^{-})\,\delta q_k+f^+_{k}(q_k, q_{k+1}, u_{k}^{+})\,\delta q_{k+1}\approx\\
  &&\approx \int_{kh}^{(k+1)h}f(q(t), \dot{q}(t), u(t))\delta{q}(t)\; dt
\end{eqnarray*}
where $\lp f^-_{k}(q_k, q_{k+1}, u_{k}^{-}), f^+_{k}(q_k, q_{k+1}, u_{k}^{+})\rp\in T_{q_k}^*Q\times T_{q_{k+1}}^*Q$ (see \cite{mawest}).

Therefore, we derive a \textbf{discrete version of the Lagrange-D'Alembert principle} given in (\ref{ldp1}):
\[
\delta \sum_{k=0}^{N-1}L_d(q_k, q_{k+1})
+\sum_{k=0}^{N-1}\left(f^{-}_{k}(q_k, q_{k+1}, u_{k}^{-})\,\delta q_k+f^+_{k}(q_k, q_{k+1}, u_{k}^{+})\,\delta q_{k+1}\right)=0,
\]
for all variations $\{\delta q_k\}_{k=0, \ldots N}$ with $\delta q_k\in T_{q_{k}}Q$ such that $\delta q_0=\delta q_N=0$.
From this principle is easy to derive the system of difference equations:
\begin{eqnarray}\label{eqw}
&&D_2L_d(q_{k-1}, q_{k})+D_1L_d(q_k, q_{k+1})\nonumber\\
&&+f^+_{k-1}(q_{k-1}, q_{k}, u_{k-1}^{+})+f^-_{k}(q_k, q_{k+1}, u_{k}^{-})=0,
\end{eqnarray}
where $k=1, \ldots, N-1$. Equations (\ref{eqw})  are called the  \textbf{forced discrete Euler-Lagrange equations} (see \cite{objuma}).

We can also approximate the cost functional (\ref{cost-1}) in a single time step $h$ by
\[
C_d(q_k, u_k^-,  q_{k+1}, u_{k}^+)
\approx \int_{kh}^{(k+1)h}C(q(t), \dot{q}(t), u(t))\; dt,
\]
yielding the \textbf{discrete cost functional}:
\[
\sum_{k=0}^{N-1}C_d(q_k, u_k^-, q_{k+1}, u_{k}^+)\; .
\]
Observe that $C_d: Q\times U\times Q\times U\rightarrow \R$.

\subsection{Fully-actuated Systems}\label{sec:fas}

In this section we assume the following condition
\begin{definition}{\rm {\bf (Fully actuated discrete system)}}
We say that the discrete mechanical control system is fully actuated if the
mappings
\begin{eqnarray*}
&f^-_{k}\big|_{(q_k, q_{k+1})}: U\to T_{q_k}^*Q, \qquad  f^-_{k}\big|_{(q_k, q_{k+1})}(u)=f^-_{k}(q_k, q_{k+1}, u),&\\
&f^{+}_{k}\big|_{(q_{k}, q_{k+1})}: U\to T_{q_{k+1}}^*Q, \qquad  f^+_{k}\big|_{(q_{k}, q_{k+1})}(u)=f^+_{k}(q_{k}, q_{k+1}, u),&
\end{eqnarray*}
are both diffeomorphisms.
\end{definition}
Define the momenta (see \cite{mawest,objuma})
\begin{eqnarray}\label{FMomenta1}
p_k&=&-D_1L_d(q_k, q_{k+1})-f^-_{k}(q_k, q_{k+1}, u_{k}^{-}),\\\label{FMomenta2}
p_{k+1}&=&D_2L_d(q_k, q_{k+1})+f^+_{k}(q_k, q_{k+1}, u_{k}^{+}).
\end{eqnarray}
Since both $f^{\pm}_{k}\big|_{(q_{k},q_{k+1})}$ are diffeomorphisms we
can express $u_{k}^{\pm}$ in terms of $(q_{k}, p_{k}, q_{k+1},
p_{k+1})$ using (\ref{FMomenta1}) and (\ref{FMomenta2}). Next, we
define a new Lagrangian ${\mathcal L}_{d}: T^*Q\times T^*Q\rightarrow
\R$ by
\vspace{0.1cm}

\begin{align}\label{eq:cLd}
  \begin{split}
    &{\mathcal L}_{d}(q_k, p_k, q_{k+1}, p_{k+1})= \\
    &=C_d(q_k\,,\,(f^{-}_{k}\big|_{(q_k,
      q_{k+1})})^{-1}(-D_1L_d-p_k)\,,\,q_{k+1}\,,\,\\
    & \hspace{50pt}(f^{+}_{k}\big|_{(q_k,
      q_{k+1})})^{-1}(-D_2L_d+p_{k+1})).
  \end{split}
\end{align}
\vspace{0.1cm}

The system is fully-actuated, consequently the Lagrangian ${\mathcal L}_{d}$ is well defined on the entire discrete space $T^*Q\times T^*Q$.

Now the discrete phase space is the Cartesian product $T^*Q\times T^*Q$ of two copies of the cotangent bundle. The definition (\ref{FMomenta1}), (\ref{FMomenta2}) gives us a matching of momenta (see \cite{mawest}) which
automatically implies
\[
D_2L_d(q_{k-1}, q_{k})+f^+_{k-1}(q_{k-1}, q_{k}, u_{k-1}^{+})=-D_1L_d(q_k, q_{k+1})-f^-_{k}(q_k, q_{k+1}, u_{k}^{-}),
\]
$k=1, \ldots, N-1$, which are the forced discrete Euler-Lagrange equations (\ref{eqw}).
In other words, the matching condition enforces that the momentum at time $k$ should be the same when evaluated
from the lower interval $[k-1, k]$ or the upper interval $[k, k+1]$. Consequently, along a solution curve there is a unique momentum at each time $t_{k}$, which can be called $p_{k}$.

The discrete Euler-Lagrange equations of motion for the Lagrangian ${\mathcal L}_{d}: T^*Q\times T^*Q\rightarrow\R$ are

\begin{eqnarray}
D_3{\mathcal L}_{d}(q_{k-1}, p_{k-1}, q_k, p_k)+D_1{\mathcal L}_{d}(q_k, p_k, q_{k+1},  p_{k+1})&=&0,\label{elca}\\
D_4{\mathcal L}_{d}(q_{k-1}, p_{k-1}, q_k, p_k)+D_2{\mathcal L}_{d}(q_k, p_k, q_{k+1}, p_{k+1})&=&0\; .\label{elcb}
\end{eqnarray}


In summary, we have obtained the discrete equations of motion for a fully-actuated mechanical optimal control problem as the discrete Euler-Lagrange equations for a Lagrangian defined on the product of two copies of the cotangent bundle. Therefore, all the preservation properties of the discrete equations (\ref{elca}) and (\ref{elcb}) are now a direct consequence of the theory of variational integrators \cite{mawest}.

\subsection{Example: optimal control problem for a mechanical Lagrangian with configuration space  $\R^{n}$}
\label{RN}

Consider the case $Q=\mathbb R^n$ and assume that $M$ is an $n\times n$ constant and
symmetric matrix. The mechanical Lagrangian
$L:\mathds{R}^{2n}\rightarrow\mathds{R}$ is defined by $L(x,\dot
x)=\frac{1}{2}\dot x^{T}M\dot x-V(x)$, where $V:\R^{n}\rightarrow\R$
is the potential function and $x\in\mathbb{R}$. The system is fully actuated and there exist no
velocity constraints. The optimal control problem is typically in
terms of boundary conditions $\lp x(0),\dot x(0)\rp$ and $\lp
x(T),\dot x(T)\rp$ for a given final time $T$.
Note that in the continuous setting we can define the momentum by the continuous Legendre transformation $\F L:TQ\Flder T^{*}Q$, $(q,\dot q)\mapsto(q,p)$:  $p=\frac{\der L}{\der\dot x}$, i.e. $p(t)=\dot x^{T}(t)\,M$. In consequence, we can define boundary constraints also in the phase space: $(x(0)\,,\,p(0)=\dot x(0)^T\,M)$ and $(x(T)\,,\,p(T)=\dot x(T)^T\,M)$.

We set the Trapezoidal discretization for the Lagrangian (see \cite{Hair}), that is, $L_{d}(x_{k},x_{k+1})=\frac{h}{2}\,L(x_{k},\frac{x_{k+1}-x_{k}}{h})+\frac{h}{2}\,L(x_{k+1},\frac{x_{k+1}-x_{k}}{h})$ where, as above, $h$ is the fixed time step and $x_{1}, x_{2}, \ldots, x_{N}$ is a sequence of elements on $\R^{n}$. Our concrete discrete Lagrangian is
\[
L_{d}(x_{k},x_{k+1})=\frac{1}{2h}(x_{k+1}-x_{k})^{T}M(x_{k+1}-x_{k})-\frac{h}{2}\lp V(x_{k})+V(x_{k+1})\rp.
\]
The control forces are $f_k^-(x_{k},x_{k+1},u_{k}^{-})\in T_{x_k}^*\R^n$ and $f_k^+(x_{k},x_{k+1},u_{k}^{+})\in T_{x_{k+1}}^*\R^n$. For sake of clarity, we are going to fix the control forces in the following manner $f^{\pm}(x_{k},x_{k+1},u_{k}^\pm)=u_{k}^{\pm}$. Looking at equations (\ref{FMomenta1}) and (\ref{FMomenta2}) is easy to obtain the associated momenta $p_{k}$ and $p_{k+1}$, namely
\begin{eqnarray*}p_k&=&\frac{1}{h}\,(x_{k+1}-x_k)^T\,M+\frac{h}{2}\,V_x(x_k)^T-u_k^{-},\\
p_{k+1}&=&\frac{1}{h}\,(x_{k+1}-x_k)^T\,M-\frac{h}{2}\,V_x(x_{k+1})^T
+ u_k^{+}.
\end{eqnarray*}
Let
$C_{d}=\frac{h}{4}\sum_{k=0}^{N-1}\left[(u_{k}^{-})^{2}+(u_{k}^{+})^{2}\right]$
be a discrete approximation of the cost function. Consequently, the
Lagrangian over $T^{*}\R^n\times T^{*}\R^n$ is
\begin{eqnarray*}
&&{\mathcal L}_{d}(x_k, p_k, x_{k+1}, p_{k+1})=\\
&&=\frac{1}{4}\sum_{k=0}^{N-1}\lp p_{k}-\lp\frac{x_{k+1}-x_{k}}{h}\rp^{T}M-\frac{h}{2}\,V_{x}(x_{k})^{T}\rp^{2}\\
&&+\frac{1}{4}\sum_{k=0}^{N-1}\lp p_{k+1}-\lp\frac{x_{k+1}-x_{k}}{h}\rp^{T}M+\frac{h}{2}\,V_{x}(x_{k+1})^{T}\rp^{2},
\end{eqnarray*}
where $V_{x}$ represents the derivative of $V$ with 	respect to the variable $x$. Applying equations (\ref{elca}) and (\ref{elcb}) to $\mathcal{L}_{d}$ we obtain the following equations:
\begin{small}
\begin{eqnarray}\label{TB1}
&&p_k-\lp\frac{x_{k+1}-x_{k-1}}{2h}\rp^T\,M=0,\\\nonumber\\\nonumber\\\label{TB2}
&&\lp p_k-(\frac{x_{k+1}-x_k}{h})^T\,M-\frac{h}{2}\,V_x(x_k)^T\rp\lp M-\frac{h^2}{2}\,V_{xx}(x_k)^T\rp\\\nonumber
&&-\lp p_k-(\frac{x_{k}-x_{k-1}}{h})^T\,M+\frac{h}{2}\,V_x(x_k)^T\rp\lp M-\frac{h^2}{2}\,V_{xx}(x_k)^T\rp=0,
\end{eqnarray}
\end{small}
where both set of equations are defined for $k=1,...,N-1$. It is quite clear that we could remove the $p_k$ dependence in equation (\ref{TB2}). However, we prefer to keep it in order to stress that the discrete variational Euler-Lagrange equations (\ref{elca}) and (\ref{elcb}) are defined in $T^*Q\times T^*Q$ ($T^*\R^n\times T^*\R^n$ in the particular case we are considering in this example).

Expressions (\ref{TB1}) and (\ref{TB2}) give $2(N-1)n$ equations for the $2(N+1)n$ unknowns $\lc x_k\rc_{k=0}^N$ , $\lc p_k\rc_{k=0}^N$. Nevertheless, the boundary conditions
\begin{eqnarray*}
&&x_0=x(0),\,\,\,\,\,p_0=p(0),\\
&&x_N=x(T),\,\,\,\,\,p_N=p(T),
\end{eqnarray*}
contribute $4n$ extra equations that convert eqs. (\ref{TB1}) and (\ref{TB2}) in a nonlinear root finding problem of $2(N-1)n$ and the same amount of unknowns.

\subsection{Underactuated Systems}\label{sec:uas}

In this section, we examine the case of underactuated systems defined
as follows:
 \begin{definition}{{\rm {\bf (Underactuated discrete system)}}}\label{UD1}
We say that the discrete mechanical control system is underactuated if the
mappings
\begin{eqnarray*}
&f^-_{k}\big|_{(q_k, q_{k+1})}: U\to T_{q_k}^*Q, \qquad  f^-_{k}\big|_{(q_k, q_{k+1})}(u)=f^-_{k}(q_k, q_{k+1}, u),&\\
&f^{+}_{k}\big|_{(q_k, q_{k+1})}: U\to T_{q_{k+1}}^*Q, \qquad  f^+_{k}\big|_{(q_k, q_{k+1})}(u)=f^+_{k}(q_k, q_{k+1}, u),&
\end{eqnarray*}
are both embeddings, that is, they are one-to-one immersions that are homeomorphisms of $U$ to its image.
\end{definition}
Under this hypothesis we deduce that ${\mathcal M}^-_{(q_k,q_{k+1})}=f^-_{k}\big|_{(q_k, q_{k+1})}(U)$,\\ ${\mathcal M}^+_{(q_k,q_{k+1})}= f^+_{k}\big|_{(q_k, q_{k+1})}(U)$ are  submanifolds of $T^*_{q_k}Q$ and $T^*_{q_{k+1}}Q$, respectively. Therefore, $f^\pm_{k}\big|_{(q_k, q_{k+1})}$ are diffeomorphisms onto its image. Moreover, $\mbox{dim}\,{\mathcal M}^-_{(q_k,q_{k+1})}=\mbox{dim}\,{\mathcal M}^+_{(q_k,q_{k+1})}=\mbox{dim}\,U$.

The set of admissible forces is restricted to the space
${\mathcal M}^-_{(q_k,q_{k+1})}\times {\mathcal
  M}^+_{(q_k,q_{k+1})}\subset T_{q_k}^*Q\times T_{q_{k+1}}^*Q$. As a
consequence, the
set of admissible momenta defined in~\eqref{FMomenta1}
and~\eqref{FMomenta2} satisfy
\begin{eqnarray*}
\lp q_k\,,\,-D_1L_d(q_k,q_{k+1})-p_k\rp&\in&\mathcal{M}^{-}_{(q_k,q_{k+1})}\subset T^*_{q_{k}}Q,\\
\lp q_{k+1}\,,\,-D_2L_d(q_k,q_{k+1})+p_{k+1}\rp&\in&\mathcal{M}^{+}_{(q_k,q_{k+1})}\subset T^*_{q_{k+1}}Q.
\end{eqnarray*}
Thus, the Lagrangian function defined in~\eqref{eq:cLd} is restricted
to these points only.  Thus, it is necessary to apply constrained
variational calculus to derive the corresponding equations (see \cite{Benito}). This is typically performed by means of constraint functions
$\Phi^-_{\alpha},\Phi^+_{\alpha}: T^*Q\times T^*Q\to \R$, $1\leq \alpha\leq n-\dim
U$. Therefore the solutions of the optimal control problem are now
viewed as the solutions of the discrete constrained problem determined
by an extended Lagrangian ${\mathcal L}_d$ and the constraints
$\Phi^{\pm}_{\alpha}$. Since $f^{\pm}\big|_{(q_k,q_{k+1})}$ are embeddings, as established in definition \eqref{UD1}, the number of constraints is determined by $n$ minus the dimension of $U$. Note that the total number of constraints, $\Phi^{\pm}_{\alpha}$, is therefore $2(n-\dim
U)$.

To solve this problem we introduce Lagrange
multipliers $(\lambda_k^-)^{\alpha}$,$(\lambda_k^+)^{\alpha}$ and
consider discrete variational calculus using the augmented Lagrangian
\begin{align*}
  \widetilde{\mathcal L}_{d}(q_k, p_k, \lambda_{k}^-, q_{k+1}, p_{k+1},
  \lambda_{k}^+) =& {\mathcal L}_{d}(q_k, p_k, q_{k+1}, p_{k+1}) \\
  & + (\lambda_k^-)^{\alpha}\Phi^-_{\alpha}(q_k, p_k, q_{k+1},
  p_{k+1}) \\
  & + (\lambda_k^+)^{\alpha}\Phi^+_{\alpha}(q_k, p_k, q_{k+1}, p_{k+1}).
\end{align*}
Observe that, in spite the constraints are functions of the Cartesian product of two copies of the cotangent bundle i.e. $\Phi_{\alpha}^{\pm}:T^*Q\times T^*Q\Flder\R$, neither $\Phi_{\alpha}^-$ depends on $p_{k+1}$ nor $\Phi_{\alpha}^+$ on $p_k$. The discrete Euler-Lagrange equations gives us the solutions of the underactuated problem.

Typically, the underactuated systems appear in an affine way that is
\begin{eqnarray*}
f^-_{k}(q_k, q_{k+1}, u_k^-)&=&A^-_{k}(q_k, q_{k+1})+B^-_{k}(q_k, q_{k+1})(u_k^-)\\
f^+_{k}(q_k, q_{k+1}, u_k^+)&=&A^+_{k}(q_k, q_{k+1})+B^+_{k}(q_k, q_{k+1})(u_k^+)
\end{eqnarray*}
where $A_k^-(q_k, q_{k+1})\in T^*_{q_k}Q$, $A_k^+(q_k, q_{k+1})\in T^*_{q_{k+1}}Q$. Moreover  $B^-_{k}(q_k, q_{k+1})\in\mbox{Lin}(U, T_{q_k}^*Q)$ and $B^+_{k}(q_k, q_{k+1})\in\mbox{Lin}(U, T_{q_{k+1}}^*Q)$ are linear maps (we assume that $U$ is a vector space and $\mbox{Lin}(E_1\,,\,E_2)$ is the set of all linear maps between $E_1$ and $E_2$). In consequence $B^-_{k}(q_k, q_{k+1})(u_k^-)\in T^{*}_{q_{k}}Q$ and $B^+_{k}(q_k, q_{k+1})(u_k^+)\in T^{*}_{q_{k+1}}Q$.

Then the constraints are deduced using the compatibility conditions:
\begin{eqnarray*}
\hbox{rank}\,B^-_{k}&=&\hbox{rank}\, \left(B^-_{k}\,;\, -D_1L_d(q_k, q_{k+1})-p_k-A^-_{k}(q_k, q_{k+1})\right),\\
\hbox{rank}\, B^+_{k}&=&\hbox{rank}\, \left(B^+_{k}\,;\, -D_2L_d(q_k, q_{k+1})+p_{k+1}-A^+_{k}(q_k, q_{k+1})\right),
\end{eqnarray*}
which imply constraints in $(q_k,q_{k+1},p_k)$ and $(q_k,q_{k+1},p_{k+1})$, respectively. The fact that $f_k^{\pm}\big|_{(q_k,q_{k+1})}$ are both embeddings implies furthermore that $\hbox{rank}\,B^-_{k}=\hbox{rank}\,B^+_{k}=\hbox{dim}\,U$.

\section{Discrete  optimal control on Lie groups}\label{RED}

An indispensable tool in the study of mechanical systems is reduction
theory. Therefore, in this work we consider its discrete analogue.
This is precisely the motivating idea of the work by Moser and
Veselov \cite{Mose}, i.e. to give a discrete analogue of
Euler-Poincar\'e reduction. The approach is to reduces the standard
second order Euler-Lagrange equations when the configuration space is
a  Lie group $G$ to first order equations on the  Lie algebra
${\mathfrak g}$.

Following the developments in $\S$ \ref{section2} assume that the
Lagrangian defined by $L_d: G\times G\rightarrow \R$ is invariant so
that
\[
L_d(g_k, g_{k+1})=L_{d}(\bar{g}g_k, \bar{g}g_{k+1})
\]
for any element $\bar{g}\in G$ and $(g_k, g_{k+1})\in G\times
G$. According to this, we can define a reduced Lagrangian $l_d:
G\rightarrow \R$ by
\[
l_d(W_k)=L_d(e, g_{k}^{-1} g_{k+1})
\]
where $W_k=g_k^{-1}g_{k+1}$ and $e$ is the identity of the Lie group $G$.

The reduced action sum is given by
\[
\begin{array}{rrcl}
S_d:& G^{N-1}&\rightarrow&\R\\
    &(W_0, \ldots, W_{N-1})&\longmapsto& \sum_{k=0}^{N-1}l_d(W_k).
\end{array}
\]
Taking variations of $S_d$ and noting that
\[
\delta W_k = -g^{-1}_k(\delta g_k) g^{-1}_k g_{k+1} + g_k^{-1}\delta
g_{k+1} = -\eta_k W_k + W_{k} \eta_{k+1},
\]
where $\eta_k = g_k^{-1}\delta g_k$,
we arrive to the discrete
Euler-Poincar\'e equations:
\[
(r_{_{W_{k}}}^*dl_d)(W_k) - (l_{_{W_{k-1}}}^*dl_d)(W_{k-1}) = 0, \ \ \
\ k = 1, ..., N-1,
\]
where $l:G\times G\Flder G$ and $r:G\times G\Flder G$ are respectively
the left and the right translations of the group (see also~\cite{BoSu}).

If we denote by $\mu_k=(r_{_{W_k}}^*dl_d)(W_k)$ then the discrete
Euler-Poincar\'e  equations are rewritten as
\begin{equation}\label{eq:dep}
  \mu_{k+1}=\mbox{Ad}^*_{_{W_{k}}}\mu_k,
\end{equation}
where $\hbox{Ad}:G\times\mathfrak{g}\Flder\mathfrak{g}$ is the adjoint
action of $G$ on $\mathfrak{g}$. Typically this equations are known as the
\emph{discrete Lie-Poisson equations } (see
\cite{BoSu,MPS1,MaPeSh}).


Consider a mechanical system  determined by a Lagrangian $l: {\mathfrak g}\rightarrow \R$, where ${\mathfrak g}$ is the Lie algebra of a Lie group $G$, which also is a $n$-dimensional vector space. The continuous external forces are defined as follows
$f: {\mathfrak g}\times U\rightarrow  {\mathfrak g}^*$. The motion of the mechanical system is described applying the following principle
\begin{equation}\label{ldp}
\delta \int^T_0 l(\xi(t))\, dt + \int^T_0 \bra f(\xi(t), u(t)),\eta(t)\ket\; dt=0,
\end{equation}
for  all variations $\delta \xi(t)$ of the form $\delta \xi(t)=\dot{\eta}(t)+[\xi(t), \eta(t)]$, where $\eta(t)$ is an arbitrary curve on the Lie algebra with $\eta(0)=0$ and $\eta(T)=0$ (see \cite{MaRa}). In addition $\bra\cdot,\cdot\ket$ is the natural pairing between $\mathfrak{g}$ and $\mathfrak{g}^{*}$.
These equations give us the controlled Euler-Poincar\'e equations:
\[
\frac{d}{dt}\left( \frac{\delta l}{\delta \xi}\right)=\hbox{ad}^*_{\xi}\left( \frac{\delta l}{\delta \xi}\right)+f,
\]
where $\hbox{ad}_{\xi}\eta=[\xi, \eta]$.

The optimal control problem consists of minimizing a given  \textbf{cost functional}:
\begin{equation}\label{cost}
\int_0^T C(\xi(t),  u(t)))\, dt,
\end{equation}
where  $C: {\mathfrak g}\times U\longrightarrow \R$.

Now, we consider the associated discrete problem. First we replace the Lie algebra ${\mathfrak g}$ by the Lie group $G$ and the continuous curves by sequences $W_0, W_1, \ldots W_N$ (since the Lie algebra is the infinitesimal version of a Lie group, its proper discretization is consequently that Lie group \cite{MaPeSh,mawest}).

The discrete Lagrangian $l_d: G\rightarrow \R$ is constructed as an approximation of the action integral, that is
\[
l_d(W_k)\approx \int_{kh}^{(k+1)h}l(\xi(t))\; dt.
\]
Let define the discrete external forces in the following way: $f^{\pm}_{k}:G\times U\Flder\mathfrak{g}^{*}$, where $U\subset\R^{m}$ for $m\leq n=\dim {\mathfrak g}$. In consequence
\[
\bra f^-_{k}(W_k,u_{k}^{-})\,,\,\eta_k\ket+\bra f^+_{k}(W_k,u_{k}^{+})\,,\,\eta_{k+1}\ket\approx \int_{kh}^{(k+1)h}\bra f(\xi(t), u(t)),\eta(t)\ket\; dt,
\]
where $(f^-_{k}(W_k,u_{k}^{-}), f^+_{k}(W_k,u_{k}^{+}))\in {\mathfrak g}^*\times {\mathfrak g}^*$ and $\eta_k\in {\frak g}$, for all $k$. In addition $\eta_0=\eta_N=0$ and $\bra\cdot,\cdot\ket$ is the natural pairing between $\al$ and $\dal$.

For sake of simplicity we are sometimes going to omit the dependence on $G\times U$ of both $f^{+}_{k}$ and $f^{-}_{k}$.

Taking all the previous into account, we derive a \textbf{discrete version of the Lagrange-D'Alembert principle for Lie groups}:
\begin{equation}\label{LDDiscreto}
\delta \sum_{k=0}^{N-1}l_d(W_k)+\sum_{k=0}^{N-1}\left(\bra f^-_{k},\eta_k\ket+\bra f^+_{k},\eta_{k+1}\ket\right)=0,
\end{equation}
for all variations $\{\delta W_k\}_{k=0, \ldots N-1}$  verifying the relation  $\delta W_k= -\eta_k  \,W_k+W_k\,\eta_{k+1}$ with $\{\eta_k\}_{k=1, \ldots N-1}$ an arbitrary sequence of elements of ${\frak g}$ which satisfies $\eta_0, \eta_N=0$ (see \cite{KM,LeMcLe2006}).

From this principle is easy to derive the system of difference equations:
\begin{align}\label{aqe}
  \begin{split}
    &l_{_{W_{k-1}}}^* d l_d(W_{k-1})-r_{_{W_{k}}}^*dl_{d}(W_{k})\\
    &+f^+_{k-1}(W_{k-1}, u_{k-1}^+)+f^-_{k}(W_k, u_k^-)=0,
  \end{split}
\end{align}
for $k=1, \ldots, N-1$, which are called the \textbf{controlled  discrete Euler-Poincar\'e equations}.

The cost functional (\ref{cost}) is approximated by
\begin{equation}\label{Cd}
  C_d(u_k^-, W_k, u^+_k)\approx \int_{kh}^{(k+1)h}C(\xi(t), u(t))\; dt,
\end{equation}
yielding the \textbf{discrete cost functional}:
\begin{equation}\label{dCostFunc}
  \mathcal{J}=\sum_{k=0}^{N-1}C_d(u_k^-, W_k , u_{k}^+)\; .
\end{equation}
Observe that now  $C_d: U\times G\times U\rightarrow \R$.

\subsection{Fully Actuated Systems}
\label{AL}
In the fully actuated case the mappings $f^{\pm}_{k}\big|_W: U\to {\mathfrak g}^*$ defined by $f^{\pm}_{k}\big|_{_{W}}(u)=f^\pm_{k}(W, u)$ are diffeomorphisms for all $W\in G$, therefore, we can construct the Lagrangian  ${\mathcal L}_d: {\mathfrak g}^*\times G\times {\mathfrak g}^*\longrightarrow \R$ by
\begin{align}  \label{Ld}
  \begin{split}
    &{\mathcal L}_d(\nu_k, W_k , \nu_{k+1})\\
    &\!\!\!\!\!\!=C_d((f^{-}_{k}\big|_{_{W_k}})^{-1}(r_{_{W_{k}}}^*dl_{d}(W_{k})\!-\!\nu_k),
    W_k, (f^{+}_{k}\big|_{_{W_k}})^{-1}(-l_{_{W_{k}}}^* dl_{d}(W_{k})+\nu_{k+1})),
  \end{split}
\end{align}
where the variables $\nu_{k},\nu_{k+1}\in\mathfrak{g}^{*}$ are defined by
\begin{align}\label{mus}
  \begin{split}
    &\nu_k=r_{W_{k}}^*dl_{d}(W_{k})-f^-_{k}(W_k, u_k^-),\\
    &\nu_{k+1}=l_{W_{k}}^* dl_{d}(W_{k})+f^+_{k}(W_k, u_k^+),
  \end{split}
\end{align}
The discrete phase space ${\mathfrak g}^*\times G\times {\mathfrak
  g}^*$ is now a mixture of two copies of the Lie algebra $\dal$ and a
Lie group $G$. This is also an example of a Lie groupoid
(\cite{groupoid}).


The discrete optimal control problem defined in
(\ref{LDDiscreto}) and (\ref{Cd}) has been reduced to a Lagrangian
one, with Lagrangian function $\mathcal{L}_{d}:{\mathfrak g}^*\times
G\times {\mathfrak g}^*\Flder\R$. In consequence,
we are able to apply discrete variational calculus to obtain
the discrete equations of motion in the phase space ${\mathfrak
  g}^*\times G\times {\mathfrak g}^*$.

Let us show how to derive these equations from a variational point of
view (see \cite{groupoid} for further details). Define first the
discrete action sum
\[
{\mathcal S}_d= \sum_{k=0}^{N-1} {\mathcal L}_d (\nu_k, W_k , \nu_{k+1}).
\]
Consider sequences of the type
$\{(\nu_k, W_k, \nu_{k+1})\}_{k=0, \ldots, N-1}$ with boundary
conditions:  $\nu_0, \nu_{N}$ and the composition
$\bar{W}=W_0\,W_1\cdot\cdot\cdot W_{N-2}\,W_{N-1}$ fixed.
Therefore an arbitrary  variation of this sequence has the form
\[
\{\nu_k(\epsilon)\,,\,h^{-1}_k(\epsilon)\,W_k\, h_{k+1}(\epsilon)\,,\,\nu_{k+1}(\epsilon)\}_{k=0, \ldots, N-1},
\]
with $\epsilon\in(-\delta, \delta)\in \R$ (both $\epsilon$ and
$\delta>0$ are real parameters) and $\nu_0(\epsilon)=\nu_0$,
$\nu_k(0)=\nu_k$ , $\nu_{N}(\epsilon)=\nu_N$, $h_k(\epsilon)\in G$ and
$h_0(\epsilon)=h_{N}(\epsilon)=e$, for all $\epsilon$. Additionally
$h_{k}(0)=e$ for all $k$.

The critical points of the discrete action sum subjected to the previous boundary conditions are characterized by
\begin{eqnarray*}
0&=&
\frac{d}{d\epsilon}\Big|_{\epsilon=0}\left(\sum_{k=0}^{N-1} {\mathcal L}_d (\nu_k(\epsilon)\,,\,h^{-1}_k(\epsilon)\,W_k\,h_{k+1}(\epsilon)\,,\,\nu_{k+1}(\epsilon))\right)\\
&=&\frac{d}{d\epsilon}\Big|_{\epsilon=0} \left\{
{\mathcal L}_d(\nu_0\,,\,W_0\,h_1(\epsilon),\,\nu_1(\epsilon))+{\mathcal L}_d(\nu_1(\epsilon)\,,\, h_1^{-1}(\epsilon)\,W_1\, h_2(\epsilon)\,,\,\nu_2(\epsilon))\right. \\
&& + \ldots + {\mathcal L}_d(\nu_{N-2}(\epsilon)\,,\,h_{N-2}^{-1}(\epsilon)\,W_{N-2}\,h_{N-1}(\epsilon)\,,\,\nu_{N-1}(\epsilon))\\
&&\left.+
{\mathcal L}_d(\nu_{N-1}(\epsilon)\,,\,h^{-1}_{N-1}(\epsilon)\,W_{N-1}\,,\,\nu_N) \right\}.
\end{eqnarray*}
Taking derivatives we obtain
  \begin{eqnarray*}
0&=&\sum_{k=1}^{N-1}\left[l_{_{W_{k-1}}}^* d{\mathcal L}_{d}\big|_{(\nu_{k-1}, \nu_{k})}(W_{k-1})-r_{_{W_{k}}}^*d{\mathcal L}_{d}\big|_{(\nu_k, \nu_{k+1})}(W_{k})\right] \delta h_k\\
&&+\sum_{k=1}^{N-1}\left[D_2{\mathcal L}_{d}\big|_{(W_{k-1})}(\nu_{k-1}, \nu_{k})+D_1{\mathcal L}_{d}\big|_{(W_{k})}(\nu_k, \nu_{k+1})\right] \delta \nu_k,
\end{eqnarray*}
where ${\mathcal L}_{d}\big|_{(W)}: {\mathfrak g}^*\times {\mathfrak g}^*\to \R$ and ${\mathcal L}_{d}\big|_{(\nu, \nu')}: G\to \R$ are defined by  ${\mathcal L}_{d}\big|_{(W)}(\nu, \nu')={\mathcal L}_{d}\big|_{(\nu, \nu')}(W)={\mathcal L}_d(\nu, W, \nu')$, where $W\in G$ and $\nu, \nu'\in {\mathfrak g}^*$.
Since $\delta h_k$ (which is defined as $\frac{d\,h_{k}}{d\epsilon}|_{\epsilon=0}$) and $\delta \nu_k$ (which is defined as $\frac{d\,\nu_{k}}{d\epsilon}|_{\epsilon=0}$), $k=1, \ldots, N-1$ are arbitrary, we deduce the following discrete equations of motion:
\begin{eqnarray}\nonumber
l_{_{W_{k-1}}}^* d{\mathcal L}_{d}\big|_{(\nu_{k-1}, \nu_{k})}(W_{k-1})-r_{_{W_{k}}}^*d{\mathcal L}_{d}\big|_{(\nu_k, \nu_{k+1})}(W_{k})&=&0,\\\label{EQUACIONES}\\\nonumber
D_2{\mathcal L}_{d}\big|_{(W_{k-1})}(\nu_{k-1}, \nu_{k})+D_1{\mathcal L}_{d}\big|_{(W_{k})}(\nu_k, \nu_{k+1})&=&0,
\end{eqnarray}
for $k=1, \ldots, N-1$. Similarly to ~\S\ref{sec:fas} we obtain the control
inputs $u_k^-$ and $u_k^+$ using~\eqref{mus}.

\subsection{Underactuated Systems}

The underactuated case can now be considered by adding of
constraints. Similarly to~\S\ref{sec:uas} underactuation restricts the
control forces to lie in a subspace spanned by vectors $\{e^s\}$ of
the basis $\{e^s, e^{\sigma}\}$ of ${\mathfrak g}^*$, where $\lc
s,\sigma\rc=1,...,n$. Then
\begin{eqnarray*}
f^-_{k}(W_k, u_k^-)&=& a^-_{k}(W_k)+(b^-_{k}(W_k,u_k^-))_s e^s,\\
f^+_{k}(W_k, u_k^+)&=& a^+_{k}(W_k)+(b^+_{k}(W_k,u_k^+))_s e^s,
\end{eqnarray*}
where $a^-_{k}(W_k), a^+_{k}(W_k)\in {\mathfrak g}^*$ and
$(b^-_{k}(W_k,u_k^-))_s, (b^+_{k}(W_k,u_k^+))_s \in \R$, for all
$s$. Additionally, the embedding condition implies that $\hbox{rank}\,b_k^-=\hbox{rank}\,b_k^+=\hbox{dim}\,U$. Then, taking the dual basis $\{e_s, e_{\sigma}\}$, we induce the
following constraints:
\begin{subequations}\label{expli}
\begin{align}
\Phi^{-}_{\sigma}(\nu_k, W_k , \nu_{k+1})&=\langle
r_{_{W_{k}}}^*dl_{d}(W_{k}) - \nu_k - a^-_{k}(W_k), e_{\sigma}\rangle=0,\\
\Phi^{+}_{\sigma}(\nu_k, W_k , \nu_{k+1})&=\langle
\nu_{k+1} - l_{_{W_{k}}}^*dl_{d}(W_{k}) - a^+_{k}(W_k), e_{\sigma}\rangle=0.
\end{align}
\end{subequations}
Observe in \eqref{expli} that, even though the constraints are functions
$\Phi_{\sigma}^{\pm}:\dal\times G\times\dal\Flder\R$, neither
$\Phi_{\sigma}^-$ depends on $\nu_{k+1}$ nor $\Phi_{\sigma}^+$ on
$\nu_k$. Once we have defined the constraints we can implement the
Lagrangian multiplier rule in order to solve the underactuated
problem. Namely, we define de extended Lagrangian as:
\begin{align}
  \begin{split}
    \tilde{\mathcal L_d}(\nu_k, \lambda_k^{-}, W_k , \nu_{k+1}, \lambda_k^{+})=& \mathcal{L}_d (\nu_k, W_k, \nu_{k+1})\\
    &+(\lambda_k^{-})^{\sigma}\Phi^{-}_{\sigma}(\nu_k, W_k ,
    \nu_{k+1})\\
    &+(\lambda_k^{+})^{\sigma}\Phi^{+}_{\sigma}(\nu_k, W_k , \nu_{k+1}).
  \end{split}
\end{align}
Defining the discrete action sum
\[
\mathcal{S}_d^{\begin{scriptsize}\mbox{under}\end{scriptsize}}=\sum_{k=0}^{N-1}\tilde
{\mathcal L}_d(\nu_k, \lambda_k^{-}, W_k , \nu_{k+1}, \lambda_k^{+}),
\]
we obtain the underactuated discrete equations of motion
\begin{align}\label{eq:udem}
  \begin{split}
    &l_{_{W_{k-1}}}^{*}\,d\mathcal{L}_d\big|_{(\nu_{k-1},\nu_k)}(W_{k-1})-r_{_{W_{k-1}}}^{*}\,d\mathcal{L}_d\big|_{(\nu_{k},\nu_{k+1})} (W_{k})\\
    &+l_{_{W_{k-1}}}^{*}\lp(\lambda_{k-1}^{-})^{\sigma}\,d\,\Phi^{-}_{\sigma}\big|_{(\nu_{k-1},\nu_k)}( W_{k-1})+(\lambda_{k-1}^{+})^{\sigma}\,d\,\Phi^{+}_{\sigma}\big|_{(\nu_{k-1},\nu_k)}( W_{k-1})\rp\\
    &-r_{_{W_{k-1}}}^{*}\lp(\lambda_{k}^{-})^{\sigma}\,d\,\Phi^{-}_{\sigma}\big|_{(\nu_{k},\nu_{k+1})}( W_{k})+(\lambda_{k}^{+})^{\sigma}\,d\,\Phi^{+}_{\sigma}\big|_{(\nu_{k},\nu_{k+1})}( W_{k})\rp=0,\\
    &D_{2}\,
    \mathcal{L}_d\big|_{(W_{k-1})}(\nu_{k-1},\nu_k) + D_{1}\, \mathcal{L}_d\big|_{(W_{k})} (\nu_{k},\nu_{k+1}) +\left[(\lambda_{k-1}^{+})^{\sigma}- (\lambda_{k}^{-})^{\sigma} \right]e_{\sigma} = 0,\\
    &\Phi^{-}_{\sigma}(\nu_k, W_k , \nu_{k+1})=0,\\
    &\Phi^{+}_{\sigma}(\nu_k, W_k , \nu_{k+1})=0,		
  \end{split}
\end{align}
\vspace{0.3cm}
where the subscripts $(W_{k-1})$, $(W_{k})$, $(\nu_{k-1},\nu_k)$,
$(\nu_{k},\nu_{k+1})$ denoted variables that are fixed.

\section{Numerical Methods for Systems on Lie Groups}
\label{DISC}

We now put the discrete optimal control equations~\eqref{EQUACIONES}
and~\eqref{eq:udem} into a form suitable for algorithmic
implementation. The numerical methods are constructed using the
following guidelines:
\begin{enumerate}
\item good approximation of the dynamics and optimality, \label{reqs:dyn}
\item avoid issues with local coordinates \label{reqs:par}
\item guarantee for numerical robustness and convergence, \label{reqs:conv}
\item numerical efficiency. \label{reqs:eff}
\end{enumerate}
The discrete mechanics approach provides an accurate approximation of
the dynamics (requirement~\ref{reqs:dyn}) through
momentum and symplectic form preservation and good energy behavior. In
addition, we will satisfy requirement~\ref{reqs:par} for systems on
Lie groups by lifting the optimization to the Lie algebra through a
retraction map that will be defined in this section. The resulting
algorithms are numerically robust in the
sense that there are no issues with coordinate singularities and the
dynamics and optimality conditions remain close to their continuous
counterparts even at big time steps. Yet, as with any other nonlinear
optimization scheme it is difficult to formally claim that the
algorithm will always
converge (requirement~\ref{reqs:conv}). Nevertheless, in practice there
are only isolated cases for underactuated systems that fail to converge. A
remedy for such cases has been suggested in~\cite{KM}. In
general, the resulting algorithms require a small number of
iterations, e.g. between 10 and 20 to converge (requirement~\ref{reqs:eff}).

The optimization variables $W_k$ are regarded as small displacements on
the Lie group. Thus, it is possible to express each term through a Lie
algebra element that can be regarded as the averaged velocity of this
displacement. This is accomplished using a {\bf retraction
  map} $\tau: {\mathfrak g}\to G$ which is an analytic local
diffeomorphism around the identity such that $\tau(\xi)\tau(-\xi)=e$,
where $\xi\in\mathfrak g$. Two standard choices for $\tau$ are
employed in this work: the exponential map, and the Cayley map.

Regarding $\xi$ as a velocity we set the discrete Lagrangian $l_d:G\to
\R$ to
\begin{eqnarray*}
  l_d(W_k)=h\,l(\xi_k),
\end{eqnarray*}
where $\xi_k=\tau^{-1}(g_k^{-1}g_{k+1})/h=\tau^{-1}(W_k)/h$. The difference
$g_{k}^{-1}\,g_{k+1}\in G$, which is an element of a nonlinear space,
can now be represented by the vector $\xi_{k}$ in order to enable
unconstrained optimization in the linear space $\mathfrak{g}$ for
optimal control purposes.

The variational principle will now be expressed in terms of the chosen
map $\tau$. The resulting discrete mechanics will thus involve the
derivatives of the map which we define next (see also~\cite{Rabee,MunteKaas,KM}):
\begin{definition}\label{Retr}
  {\rm Given a map $\tau:\mathfrak{g}\Flder G$, its {\bf right trivialized tangent} $\mbox{d}\tau_{\xi}:\mathfrak{g}\Flder\mathfrak{g}$ and is {\bf inverse} $\mbox{d}\tau_{\xi}^{-1}:\mathfrak{g}\Flder\mathfrak{g}$, are such that for $g=\tau(\xi)\in G$ and $\eta\in\mathfrak{g}$, the following holds}
  \begin{eqnarray*}
    &&\der_{\xi}\tau(\xi)\,\eta=\mbox{d}\tau_{\xi}\,\eta\,\tau(\xi),\\
    &&\der_{\xi}\tau^{-1}(g)\,\eta=\mbox{d}\tau^{-1}_{\xi}(\eta\,\tau(-\xi)).
  \end{eqnarray*}
\end{definition}
Using these definitions, variations $\delta\xi$ and $\delta g$ are
constrained by
\[
\delta \xi_k=\hbox{d}\tau^{-1}_{h\xi_k}(-\eta_k+\hbox{Ad}_{\tau(h\xi_k)}
\eta_{k+1})/h,
\]
where $\eta_k=g_k^{-1}\delta g_k$, which is obtained by
straightforward differentiation of $\xi_k=\tau^{-1}(g_{k}^{-1}\,g_{k+1})/h$.

The retraction map $\tau $ choices are:
\vspace{0.2cm}

a) The exponential map $\e:\al\Flder G$, defined by $\exp(\xi)=\gamma(1)$, with $\gamma:\R\Flder G$ in the integral curve through the identity of the vector field associated with $\xi\in\al$ (hence, with $\dot\gamma(0)=\xi$). The right trivialized derivative and its inverse are defined by
\begin{eqnarray*}
\mbox{d}\e_{x}\,y&=&\sum_{j=0}^{\infty}\frac{1}{(j+1)!}\,\ad_{x}^{j}\, y,\\
\mbox{d}\e_{x}^{-1}\,y&=&\sum_{j=0}^{\infty}\frac{B_{j}}{j!}\,\ad_{x}^{j}\, y,
\end{eqnarray*}
where $B_{j}$ are the Bernoulli numbers (see \cite{Hair}). Typically, these expressions are truncated in order to achieve a desired order of accuracy.
\vspace{0.2cm}

b) The Cayley map $\ca:\al\Flder G$ is defined by $\ca(\xi)=(e-\frac{\xi}{2})^{-1}(e+\frac{\xi}{2})$ and is valid for a general class of quadratic groups (see \cite{Hair}) that include the groups of interest in this paper (e.g. $SO(3)$, $SE(2)$ and $SE(3)$). Its right trivialized derivative and inverse are defined by
\begin{eqnarray*}
\mbox{d}\ca_{x}\,y&=&(e-\frac{x}{2})^{-1}\,y\,(e+\frac{x}{2})^{-1},\\
\mbox{d}\ca_{x}^{-1}\,y&=&(e-\frac{x}{2})\,y\,(e+\frac{x}{2}).
\end{eqnarray*}

Next, the discrete forces and cost function are defined through a
trapezoidal approximation, i.e.
\[
f^{\pm}_{k}(\xi_k, u_k^{\pm})=\frac{h}{2}\,f(\xi_k,u_k^{\pm}),
\]
and
\[
C_d(u_k^-,\xi_k, u_k^+)=\frac{h}{2}\,C(\xi_k,u_k^-)+\frac{h}{2}\,C(\xi_k,u_k^+),
\]
respectively. With the choice of a retraction map and the trapezoidal
rule the equations of motion~\eqref{eq:dep} become
\begin{eqnarray*}
  &&\mu_k-\mbox{Ad}^*_{\tau(h\xi_{k-1})}\mu_{k-1}=\frac{h}{2}\,f(\xi_k, u_k^-)+\frac{h}{2}\,f(\xi_{k-1}, u_{k-1}^+),\\
  &&\mu_k=(\mbox{d}\tau^{-1}_{h\xi_k})^*\partial_{\xi}l(\xi_k),\\
  &&g_{k+1}=g_k\tau(h\xi_k),
\end{eqnarray*}
while the momenta defined in~\eqref{mus} take the form
\begin{align}
  &\nu_k=\mu_k-\frac{h}{2}\,f(\xi_k, u_k^-),\label{nu1}\\
  &\nu_{k+1}=\hbox{Ad}^*_{\tau(h\xi_k)}\mu_k+\frac{h}{2}\,f(\xi_k, u_k^+).\label{nu2}
\end{align}
Finally, define the Lagrangian ${\ell}_{d}:\dal\times\al\times\dal\to \R$ such that
\[
\ell_{d}(\nu,\xi,\nu') = \mathcal{L}_{d}(\nu, \tau(h\xi), \nu').
\]
Note that the Lagrangian is well-defined only on $\mathfrak g^*\times
\mathfrak U \times \mathfrak g^*$, where $\mathfrak U\subset\mathfrak
g$ is an open neighborhood around the identity for which $\tau$ is a
diffeomorphism. To make the notation as simple as possible we retain the
Lagrangian definition to the full space $\mathfrak g^*\times
\mathfrak g \times \mathfrak g^*$.

The optimality conditions corresponding to~\eqref{EQUACIONES} become
\begin{eqnarray}\label{EQUACIONES2Pr}
  (\mbox{d}\tau^{-1}_{-h\xi_{k-1}})^{*}\,d\,\ell_{d}\big|_{(\nu_{k-1},\nu_{k})}(\xi_{k-1})-(\mbox{d}\tau^{-1}_{h\xi_{k}})^{*}\,d\,\ell_d\big|_{(\nu_{k},\nu_{k+1})}(\xi_{k})&=&0,\\\label{EQUACIONES2}
D_{2}\,\ell_{d}\big|_{(\xi_{k-1})}(\nu_{k-1},\nu_{k})+D_{1}\,\ell_{d}\big|_{(\xi_{k})}(\nu_{k},\nu_{k+1})&=&0,
\end{eqnarray}
for $k=0,...,N-1$. Here,
$\ell_{d}\big|_{(\xi)}(\nu,\nu^{\prime})=\ell_{d}\big|_{(\nu,\nu^{\prime})}(\xi)=\ell_{d}(\nu,\xi,\nu^{\prime}).$
Equations (\ref{EQUACIONES2Pr}) and (\ref{EQUACIONES2}) can be also
obtained from (\ref{EQUACIONES}) employing Lemma~\ref{lemaPr} and
Lemma~\ref{lemaPr2} in Appendix A.

In the underactuated case we define
\begin{align}
  \begin{split}
    \tilde\ell_{d}(\nu,\xi,\nu',\lambda^-,\lambda^+) =&
    \mathcal{L}_{d}(\nu, \tau(h\xi), \nu')\\
    & +
    (\lambda^{-})^{\sigma}\Phi^{-}_{\sigma}\big|_{(\nu,\nu')}(
    \tau(h\xi)) +(\lambda^{+})^{\sigma}\Phi^{+}_{\sigma}\big|_{(\nu,\nu')}(\tau(h\xi)),
  \end{split}
\end{align}
and from~\eqref{eq:udem} obtain the equations
\begin{align}\label{eq:lieopt_uc}
  \begin{split}
    &
    (\mbox{d}\tau^{-1}_{-h\xi_{k-1}})^{*}\,d\,\tilde\ell_{d}\big|_{(\nu_{k-1},\nu_{k},\lambda_{k-1}^\pm)}(\xi_{k-1})-(\mbox{d}\tau^{-1}_{h\xi_{k}})^{*}\,d\,\tilde\ell_d\big|_{(\nu_{k},\nu_{k+1},\lambda_k^{\pm})}(\xi_{k})
    = 0,\\
    &D_{2}\,
    \mathcal{L}_d\big|_{\tau(h\xi_{k-1})}(\nu_{k-1},\nu_k) + D_{1}\,
    \mathcal{L}_d\big|_{\tau(h\xi_k)} (\nu_{k},\nu_{k+1})
    +\lambda_{k-1}^{+} - \lambda_{k}^-  = 0,\\
    &\Phi^{-}_{\sigma}(\nu_k, \tau(h\xi_k) , \nu_{k+1})=0,\\
    &\Phi^{+}_{\sigma}(\nu_k, \tau(h\xi_k) , \nu_{k+1})=0,		
  \end{split}
\end{align}
where we employed the notation $\lambda^\pm := (\lambda^\pm)^\sigma
e_{\sigma}$.


\vspace{0.3cm}

{\it Boundary Conditions.}: Establishing the exact relationship
between the discrete and continuous momenta, $\mu_{k}$ and
$\mu(t)=\der_{\xi}l(\xi(t))$, respectively, is particularly important
for properly enforcing boundary conditions that are given in terms of
continuous quantities. The following equations relate the momenta at
the initial and final times $t=0$ and $t=T$ and are used to transform
between the continuous and discrete representations:
\begin{eqnarray*}
\mu_{0}-\der_{\xi}l(\xi(0))&=&\frac{h}{2}\,f(\xi(0),u_{0}^{-}),\\
\der_{\xi}l(\xi(T))-\mbox{Ad}_{\tau(h\xi_{N-1})}^{*}\,\mu_{N-1}&=&\frac{h}{2}\,f(\xi(T),u_{N}^{+}).
\end{eqnarray*}
which also corresponds to the relations  $\nu_{0}=\der_{\xi}l(\xi(0))$
and $\nu_{N}=\der_{\xi}l(\xi(T))$. These equations can also be
regarded as structure-preserving {\bf velocity boundary conditions},
i.e., for given fixed velocities $\xi(0)$ and $\xi(T)$.

The exact form of the previous equations depends on the choice of
$\tau$. This choice will also influence the computational efficiency of the
optimization framework when the above equalities are enforced as
constraints. The numerical procedure to compute the trajectory is
summarized as follows:

\begin{algorithm}{\bf Optimal control}\label{algo:oc}
\begin{enumerate}
  \setlength{\itemsep}{-0pt}
\item[] Data: group $G$; mechanical Lagrangian $l$; control
  functions $a$, $b$; cost function $C$; final time $T$; number of
  segments $N$.
\item \ \ \ Input: boundary conditions $(g(0),\xi(0))$ and
  $(g(T),\xi(T))$.
\item \ \ \ Set momenta $\nu_0=\der_{\xi}l(\xi(0))$ and $\nu_N=\der_{\xi}l(\xi(T))$
\item \ \ \ Solve for
  $(\xi_0,...,\xi_{N-1},\nu_1,...,\nu_{N-1},\lambda_1^{\pm},...,\lambda_{N-1}^{\pm})$ the relations:\\
  $\left\{\begin{array}{l}
      \text{equations}~\eqref{eq:lieopt_uc} \text{ for all } k=1,...,N-1, \\
      \tau^{-1}\lp\tau(h\xi_{N-1})^{-1}...\tau(h\xi_{0})^{-1}\,g(0)^{-1}g(T)\rp=0
    \end{array}\right.$
  \item \ \ \ Output: optimal sequence of velocities
  $\xi_0,...,\xi_{N-1}$.
\item \ \ \ Reconstruct path $g_0,...,g_N$ by $g_{k+1} =
  g_k\tau(h\xi_k)$ for $k=0,...,N-1$.
\end{enumerate}
\end{algorithm}
The solution is computed using root-finding procedure such as Newton's
method. If the initial guess does not satisfy the dynamics we
recommend to use a Levenberg-Marquardt algorithm which has slower but
more robust convergence.

\subsection{Example: optimal control effort}
\label{Example}
Consider a Lagrangian consisting of the kinetic energy only
\[
l(\xi)=\frac{1}{2}\bra\I(\xi)\,,\,\xi\ket,
\]
full unconstrained actuation, no potential or external forces and
no velocity constraint. The map $\I:\al\Flder\dal$ is called the
inertia tensor and is assumed full rank.

In the fully actuated case we have $f(\xi_{k},u^{\pm}_{k})\equiv u^{\pm}_{k}$.
We consider a minimum effort control problem, i.e.
\[
C(\xi,u) = \frac{1}{2}\|u\|^2.
\]
The optimal control problem for fixed initial and final states
$(g(0)\,,\,\xi(0))$ and $(g(T)\,,\,\xi(T))$ can now be summarized as:
\vspace{0.4cm}

{\bf Compute:} $\xi_{0:N-1}$ , $u^{\pm}_{0:N}$,
\vspace{0.4cm}

{\bf minimizing:} $\frac{h}{4}\sum_{k=0}^{N-1}\lp\parallel u_{k}^{-}\parallel^{2}+\parallel u_{k}^{+}\parallel^{2}\rp,$
\vspace{0.4cm}

{\bf subject to:}
\vspace{0.2cm}

$\mu_{0}-\I(\xi(0))=\frac{h}{2}\,u_{0}^{-}$,
\vspace{0.1cm}

$\mu_{k}-\Ad^{*}_{\tau(h\xi_{k-1})}\,\mu_{k-1}=h(u_{k}^{-}+u_{k-1}^{+}),\hspace{2.5cm}
k=1,...,N-1,$
\vspace{0.1cm}

$\I(\xi(T))-\Ad^{*}_{\tau(h\xi_{N-1})}\,\mu_{N-1}=\frac{h}{2}\,u_{N}^{+}$,
\vspace{0.1cm}

$\mu_{k}=\dmuno\,\I(\xi_{k})$,
\vspace{0.1cm}

$g_{k+1}=g_{k}\,\tau(h\xi_{k}),\hspace{6cm} k=0,...N-1,$
\vspace{0.1cm}

$\tau^{-1}(g_{N}^{-1}\,g(T))=0.$
\vspace{0.2cm}

The optimality conditions for this problem are derived as
follows. The Lagrangian becomes
\begin{small}
  \[
  \ell_{d}(\nu_{k},\xi_{k},\nu_{k+1})=\frac{1}{4h}\sum_{k=0}^{N-1}\lp\parallel\nu_{k}-(\mbox{d}\,\tau^{-1}_{h\xi_{k}})^{*}\I(\xi_{k})\parallel^{2}+\parallel\nu_{k+1}-(\mbox{d}\,\tau^{-1}_{-h\xi_{k}})^{*}\I(\xi_{k})\parallel^{2}\rp,
  \]
\end{small}
where the momentum has been computed according to
\begin{equation}\label{NU}
  \nu_{k}=\frac{1}{2}\lp(\mbox{d}\,\tau^{-1}_{h\xi_{k}})^{*}\I(\xi_{k})+(\mbox{d}\,\tau^{-1}_{-h\xi_{k-1}})^{*}\I(\xi_{k-1})\rp,
\end{equation}
Thus the optimality conditions become
\begin{eqnarray*}
  &&(\mbox{d}\,\tau^{-1}_{h\xi_{k}})^{*}\,d\ell_{d}\big|_{(\nu_{k},\nu_{k+1})}(\xi_{k})-(\mbox{d}\,\tau^{-1}_{-h\xi_{k-1}})^{*}\,d\ell_{d}\big|_{(\nu_{k-1},\nu_{k})}(\xi_{k-1})=0,\\
  &&k=1,...,N-1,\\\\
  &&\tau^{-1}\lp\tau(h\xi_{N-1})^{-1}...\tau(h\xi_{0})^{-1}\,g_{0}^{-1}g(T)\rp=0.
\end{eqnarray*}
It is important to note that these last two equations define $N\cdot n$
equations in the $N\dot n$ unknowns $\xi_{0:N-1}$. A solution can be
found using nonlinear root finding. Once $\xi_{0:N}$ have been
computed, is possible to obtain the final configuration $g_{N}$ by
reconstructing the curve by these velocities. Beside, the boundary
condition $g(T)$ is enforced through the relation
$\tau^{-1}(g_{N}^{-1}\,g(T))=0$ without the need to optimize over any
of the configurations $g_{k}$.

\subsection{Extension: the configuration-dependent case}\label{Potencial}
The developed framework can be extended to a
configuration-dependent Lagrangian
$L:G\times\al\Flder\R$, for instance defined in terms of a kinetic
energy $K:\mathfrak g\rightarrow\mathbb R$ and potential energy
$V:G\Flder\R$ according to
\[
L(g, \xi)=K(\xi)-V(g),
\]
where $g\in G$ and $\xi\in\al$. The controlled Euler-Poincar\'e
equations are in this case
\begin{eqnarray*}
  &&\dot\mu-\mbox{ad}_{\xi}^{*}\mu=-l_g^{*}\,\der_{g}\,V(g)+f,\\
  &&\mu=\der_{\xi}K(\xi),\\
  &&\dot g=g\,\xi,
\end{eqnarray*}
where the external forces are defined as $f:G\times\al\times U\Flder\dal$. Our discretization choice $L_{d}:G\times G\Flder\R$ will be (recall that $\xi_k=\tau^{-1}(g_k^{-1}g_{k+1})/h$)
\begin{eqnarray*}
L_d(g_k, g_{k+1})&=&\frac{h}{2}L(g_k, \xi_k)+\frac{h}{2}L(g_{k+1}, \xi_k)\\
&=& h\,K(\xi_k)-h\,\frac{V(g_k)+V(g_{k+1})}{2},
\end{eqnarray*}
while the $G$-dependent discrete forces now become
\[
f_k^{-}(g_k,\xi_k,u_k^-)=\frac{h}{2}\,f(g_k,\xi_k,u_k^-), \ \ \ \ f_k^{+}(g_{k+1},\xi_k,u_k^+)=\frac{h}{2}\,f(g_{k+1},\xi_k,u_k^+).
\]
This leads to the discrete equations
\begin{eqnarray*}
&&\mu_k-\mbox{Ad}^*_{\tau(h\xi_{k-1})}\mu_{k-1}=-h\,l_{g_k}^*\partial_gV(g_k)
\\&& + \frac{h}{2}\,f(g_k, \xi_k, u_k^-)+\frac{h}{2}\,f(g_{k}, \xi_{k-1}, u_{k-1}^+),\\\\
&&\mu_k=(\hbox{d}\tau^{-1}_{h\xi_k})^*\partial_{\xi}K(\xi_k),\\
&&g_{k+1}=g_k\tau(h\xi_k).
\end{eqnarray*}

The  momenta become
\begin{eqnarray*}
\nu_k&=&\mu_k+\frac{h}{2}\,l_{g_k}^*\,\der_g\,V(g_k)-\frac{h}{2}\,f(g_k, \xi_k, u_k^-),\\
\nu_{k+1}&=&\hbox{Ad}^*_{\tau(h\xi_k)}\mu_k-\frac{h}{2}\,l_{g_{k+1}}^*\,\partial_g\,V(g_{k+1})+\frac{h}{2}\,f(g_{k+1},  \xi_k, u_k^+).
\end{eqnarray*}
In consequence, we can define a discrete Lagrangian
\[
\mathfrak{ L}_{d}:\dal\times\ G\times\al\times\dal\to \R,
\]
depending on the variables $(\nu_k, g_{k},\xi_{k}, \nu_{k+1})$ which discrete equations of motion will be a mixture between (\ref{EQUACIONES}) and (\ref{EQUACIONES2Pr}), (\ref{EQUACIONES2}), namely

\begin{align*}
  & D_{2}\,\mathfrak{ L}_{d}\big|_{(g_{k-1},\xi_{k-1})}(\nu_{k-1},\nu_{k})+D_{1}\,\mathfrak{ L}_{d}\big|_{(g_{k},\xi_{k})}(\nu_{k},\nu_{k+1})=0,\\
  & \lp l_{_{g_{k-1}}}^* d\,\mathfrak{ L}_{d}\big|_{(\nu_{k-1},\xi_{k-1},\nu_{k})}(g_{k-1})+r_{_{g_{k}}}^*d\,\mathfrak{ L}_{d}\big|_{(\nu_k,\xi_{k},\nu_{k+1})}(g_{k})\rp\\
  &+\lp(\mbox{d}\tau^{-1}_{-h\xi_{k-1}})^{*}\,d\,\mathfrak{ L}_{d}\big|_{(\nu_{k-1},g_{k-1},\nu_{k})}(\xi_{k-1})-(\mbox{d}\tau^{-1}_{h\xi_{k}})^{*}\,d\,\mathfrak{ L}_{d}\big|_{(\nu_{k},g_{k},\nu_{k+1})}(\xi_{k})\rp=0.
\end{align*}

\section{Applications}\label{Vehiculo}
\subsection{Underwater Vehicle}\label{sec:uv}
\begin{figure}[h]
  \includegraphics[width=5in]{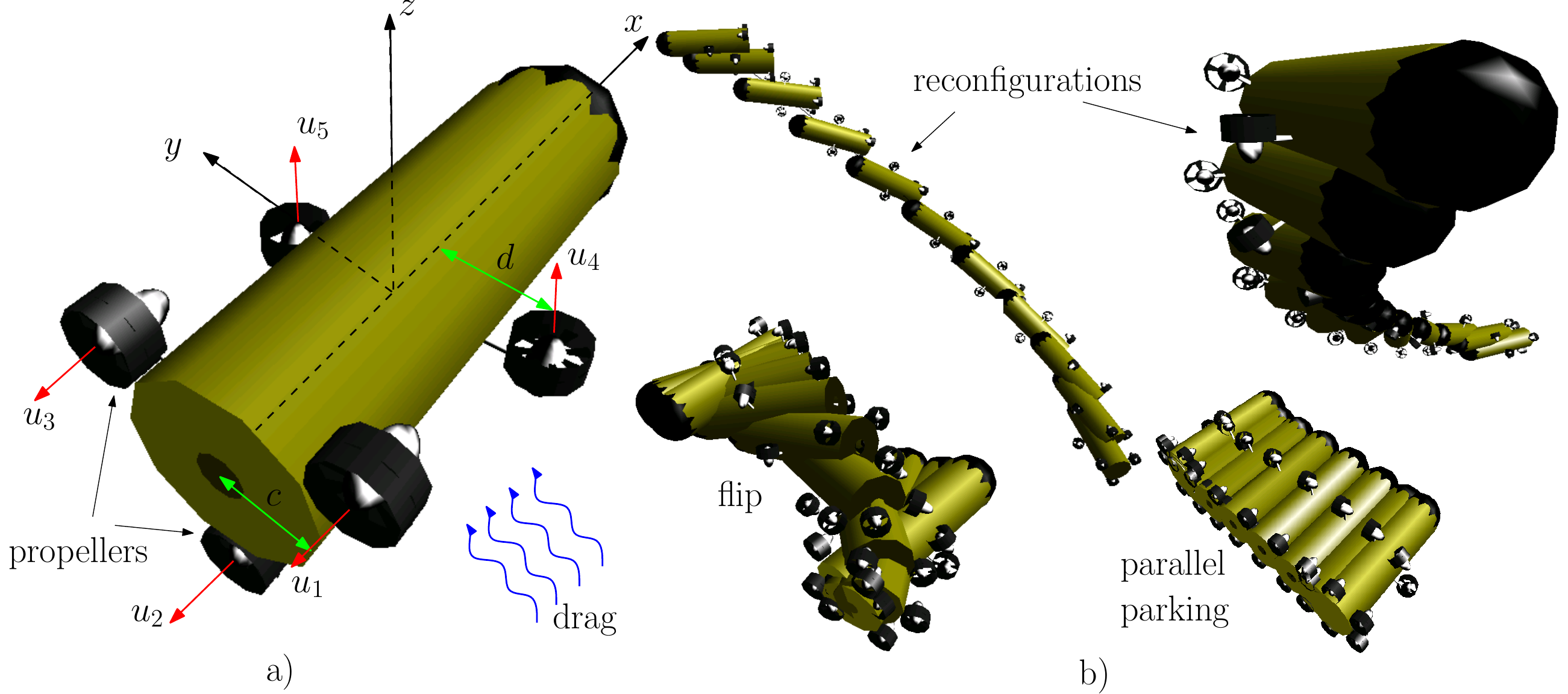}
  \caption{\footnotesize An underwater vehicle model (a) and a various
    computed optimal trajectories between chosen states (b).  Only a
    few frames along the path are shown for clarity.}\label{fig:uuv}
\end{figure}

We illustrate the developed algorithm with an application to a
simulated unmanned underwater vehicle. Figure \eqref{fig:uuv} shows
the model equipped with five thrusters which produce forces and torques in
all directions but the body-fixed ``y''-axis. Since the input
directions span only a five-dimensional subspace the problem is solved
through the underactuated framework.

\begin{figure}[h]\label{fig:uuv_plots}
  \begin{tabular}{@{}c@{}@{}c@{}@{}c@{}}
    \includegraphics[width=1.65in]{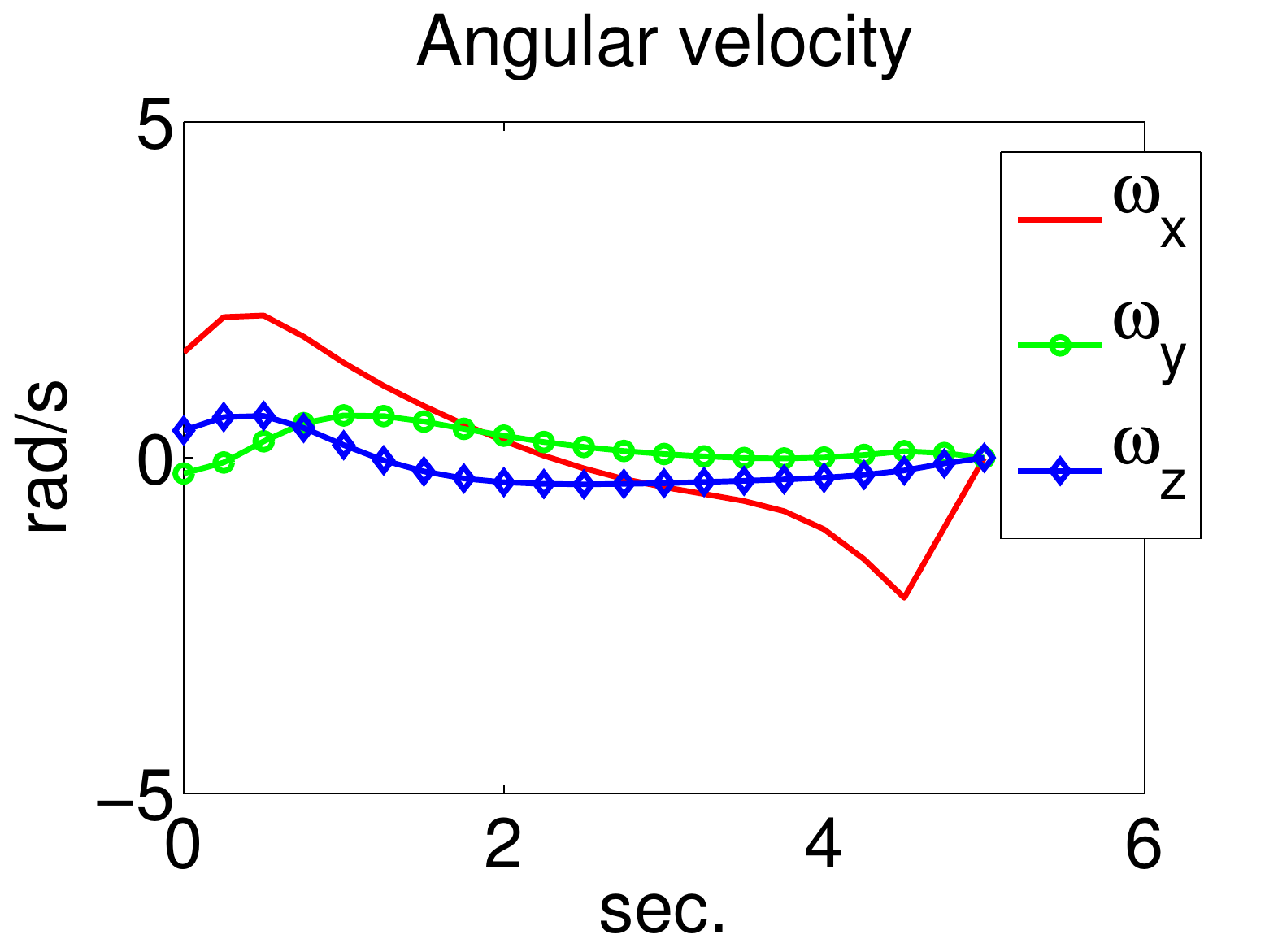}
    & \!\!\!\includegraphics[width=1.65in]{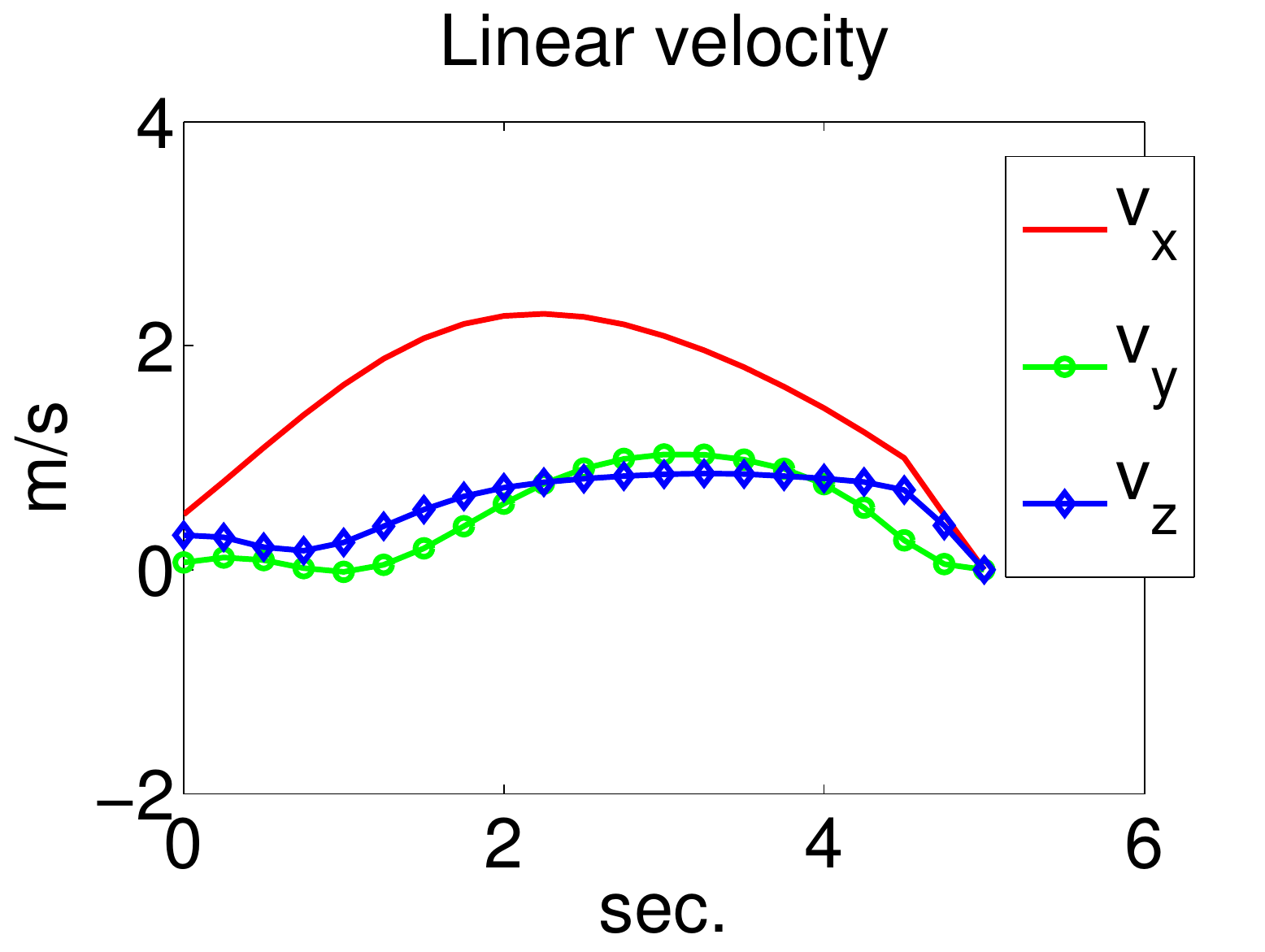} & \!\includegraphics[width=1.65in]{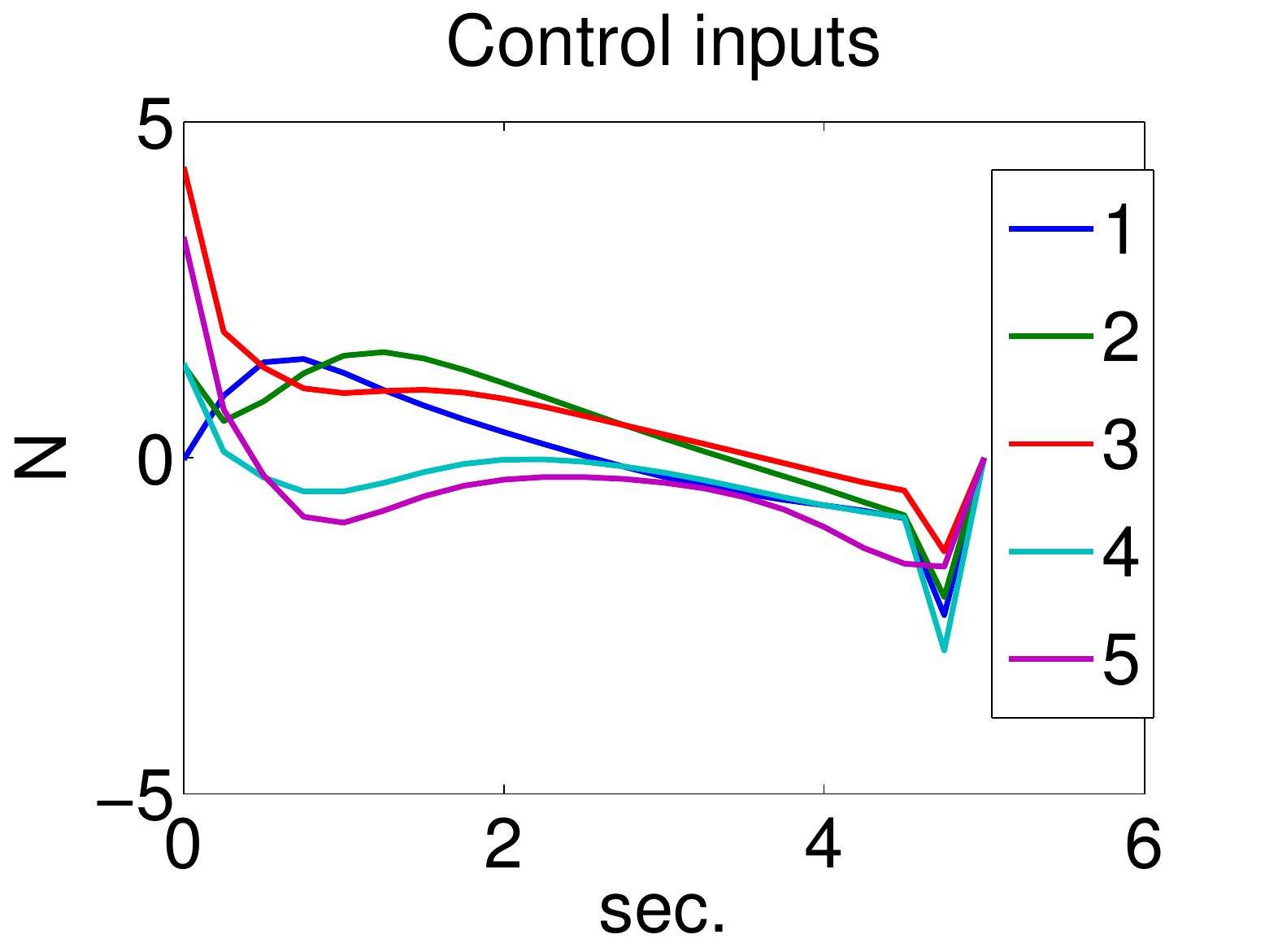}
  \end{tabular}
  \caption{\footnotesize Details of the computed optimal path for the
    reconfiguration maneuver given in
    Figure \eqref{fig:uuv}.}\label{fig:plots}
\end{figure}

The vehicle configuration space is $G=SE(3)$. We make the
identification $SE(3)\sim SO(3)\times\R^{3}$ using elements $R\in
SO(3)$ and $x\in\R^{3}$ through
\[
g=\lp\begin{array}{cc}R&x\\
0_{3\times 3}&1
\end{array}
\rp,\,\,\,\,g^{-1}=\lp\begin{array}{cc}
R^{T}&-R^{T}x\\
0_{3\times 3}&1
\end{array}
\rp,
\]
where $g\in SE(3)$. Elements of the Lie algebra $\xi\in\setr$ are
identified with body-fixed angular and linear velocities denoted
$\omega\in\R^{3}$ and $v\in\R^{3}$, respectively, through
\[
\xi=\lp\begin{array}{cc}
\hat\omega&v\\
0_{3\times 3}&0
\end{array}
\rp,
\]
where the map $\hat\cdot:\R^{3}\Flder\alg$ is defined by
\begin{equation}\label{hatso}
  \hat\omega=\lp\begin{array}{ccc}
    0&-\omega_{3}&\omega_{2}\\
    \omega_{3}&0&-\omega_{1}\\
    -\omega_{2}&\omega_{1}&0
  \end{array}
  \rp.
\end{equation}
The algorithm is thus implemented in terms of vectors in
$\mathbb{R}^6$ rather than matrices in $\mathfrak {se} (3)$.

The map $\tau=\ca:\mathfrak {se} (3)\rightarrow SE(3)$ is chosen, instead of the exponential, since it results in more computationally
efficient implementation. It is defined by
\[
\ca(\xi)=\lp\begin{array}{cc}
\ca(\hat\omega)&\mbox{d}\ca_{\omega}\,v\\
0&1
\end{array}
\rp,
\]
where $\ca:\alg\Flder SO(3)$ is given~\footnote{note that
  $\ca$ denotes a map to either $SO(3)$ or $SE(3)$ which should be clear
  from its argument.} by
\begin{equation}\label{caySO}
  \ca(\hat\omega)=\mathbf I_{3}+\frac{4}{4+\parallel\omega\parallel^{2}}\lp\hat\omega+\frac{\hat\omega^{2}}{2}\rp,
\end{equation}
where $\mathbf I_n$ is the $n\times n$ identity matrix and
$\mbox{d}\ca:\R^{3}\Flder\R^{3}$ is defined by
\begin{equation}\label{Dtau}
  \mbox{d}\ca_{\omega}=\frac{2}{4+\parallel\omega\parallel^{2}}(2\mathbf
  I_3+\hat\omega).
\end{equation}
The matrix representation of the right-trivialized tangent inverse
$\operatorname{d\tau}^{-1}_{(\omega,v)}:\mathbb{R}^3\times\mathbb{R}^3\rightarrow
\mathbb{R}^3\times\mathbb{R}^3$ becomes
\begin{align}
  &[\operatorname{dcay}^{-1}_{(\omega,v)}] = \left[ \begin{array}{cc}
      \mathbf{I}_3 - \frac{1}{2}\widehat\omega +
      \frac{1}{4}\omega\omega^T & \mathbf{0}_3 \\
      -\frac{1}{2}\left(\mathbf{I}_3 -
        \frac{1}{2}\widehat\omega\right)\widehat{v} & \mathbf{I}_3 -
      \frac{1}{2}\widehat\omega \end{array}
  \right]. \label{eq:se3_dcayinv}
\end{align}

The vehicle inertia tensor $\mathbb{I}$ is computed assuming
cylindrical mass distribution with mass $m=3$kg. The control basis
vectors are $\{e_s\}_{s=1}^{5}=\{\mathbf e_1,\mathbf e_2,
\mathbf e_3, \mathbf e_4, \mathbf e_5\}$, while the non-actuated
direction is $e_\sigma=\mathbf e_6$, where $\mathbf e_i$ is the
$i$-th standard basis vector of $\mathbb R^6$.
The control functions take the form
\begin{align*}
  &b(W,u)_1 = d(u_5 - u_4) , \\
  &b(W,u)_2 = c((u_1 + u_2)/2 - u_3), \\
  &b(W,u)_3 = (c\sin\frac{\pi}{3})(u_2 - u_1) , \\
  &b(W,u)_4 = u_1 + u_2 + u_3,  \\
  &b(W,u)_5 = u_4 + u_5,  \\
  &a(W) = H\tau^{-1}(W),
\end{align*}
here $H$ is a negative definite viscous drag matrix and the constants
$c,d$ are the lengths of the thrusting torque moment arms (see
Figure~\ref{fig:uuv}).

We are interested in computing a minimum control effort trajectory
between two given boundary states, i.e. conditions on both the
configurations and velocities. Such a cost function is defined
in~\S\ref{Example}. The optimal control problem is solved
using equations~\eqref{eq:lieopt_uc}. The computation is performed
using Algorithm~\ref{algo:oc}. Figure~\ref{fig:plots} shows the
computed velocities and controls for the ``reconfiguration''
trajectory shown in Figure~\ref{fig:uuv}. The algorithms requires
between 10-20 iterations depending on the boundary conditions and when
applied to $N=32$ segments.

\subsection{Discontinuous Control}\label{sec:dc}
\begin{figure}[h]
  \footnotesize
  \begin{tabular}{cc}
    \!\!\!\!\!\!\includegraphics[width=3in]{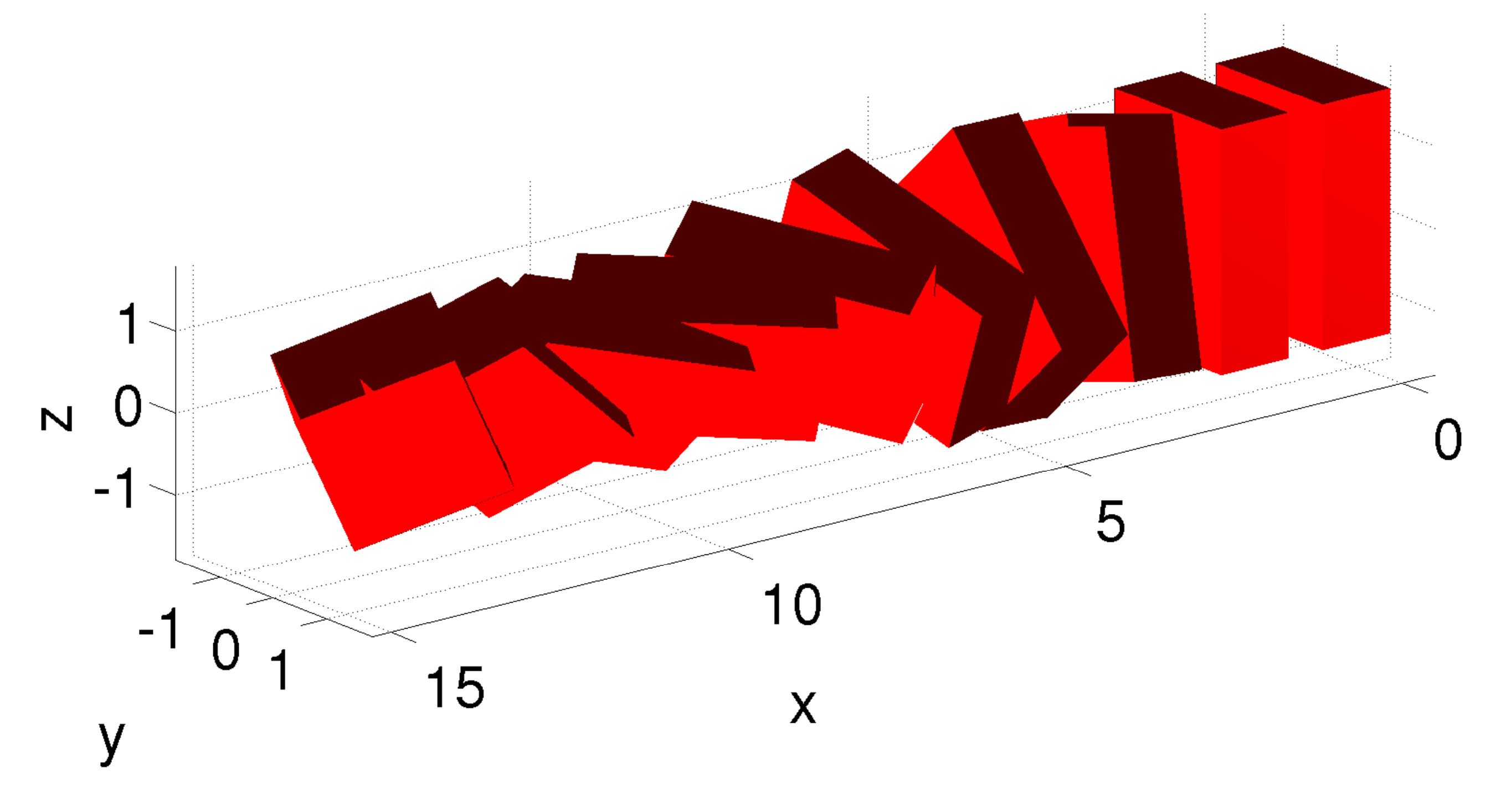}
    & \!\!\!\!\!\includegraphics[width=2in]{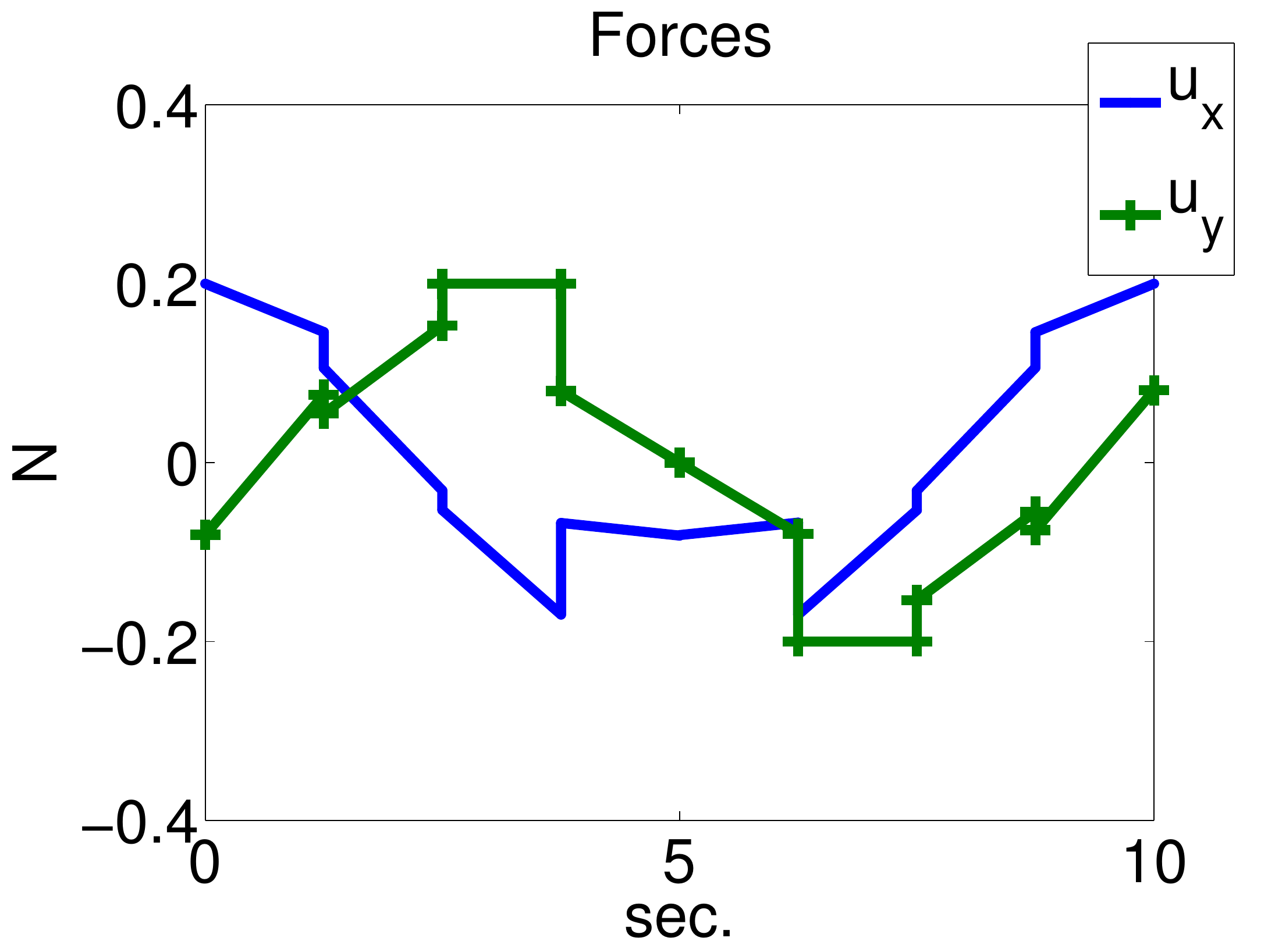} \\
    \footnotesize a) &\footnotesize b)
  \end{tabular}

  \caption{\footnotesize An optimal trajectory of an underactuated rigid body on
    $SO(3)$ (a). The body is controlled using two force inputs around
    the body-fixed $x$ and $y$ axes. An $L_1$-control cost function
    results in a discontinuous optimal trajectory (b) which our
    algorithm can handle.}\label{fig:rb}
\end{figure}

One of the advantages of employing the discrete variational framework
is the treatment of discontinuous control inputs as illustrated
in~\S\ref{TQ}. The nature of the control curve depends on the cost
function. In the standard squared control effort case (i.e. $L_2$
control curve norm employed in~\S\ref{sec:uv}) the resulting control is
smooth. Another cost
function of interest is $\int_0^T \|u(t)\|dt$ (i.e. the $L_1$ control
curve norm) which is typically imposed along with the constraints $u_{min}
\leq u(t) \leq u_{max}$. This case results in a discontinuous optimal
control curve. Our formulation can handle such problems easily since
the terms $u_k^-$ and $u_k^+$ are regarded as the forces before and
after time $t_k$, respectively. A computed scenario of a rigid body
actuated with two control torques around its principles axes of
inertia (Fig.~\ref{fig:rb}) illustrates the discontinuous case.

\section{Extensions}\label{sec:ext}
\label{Ati}
The methods developed in the previous sections are easily adapted to other cases which are of great interest in real applications. In particular, this section will be devoted to the discussion of two important extensions: the case of optimal control problems for Lagrangians of the type $l: TM\times {\mathfrak g}\to \R$ (that is, reduction by symmetries on a trivial principal fiber bundle) and the case of nonholonomic systems. Here, $M$ denotes a smooth manifold.
Observe that the phase space $TM\times {\mathfrak g}$ unifies the
previously studied cases of a tangent bundle and a Lie algebra.

The notion of principal fiber bundle is present
in many locomotion and robotic systems~\cite{BlKrMaMu1996,bullolewis, MaOs1998}. When the configuration manifold is $Q=M\times G$, there exists a canonical splitting between variables describing the position and variables describing the
orientation of the mechanical system. Then, we distinguish the pose coordinates $g\in G$ (the elements in the Lie algebra will be denoted by $\xi\in\al$),  and the variables describing the
internal shape of the system, that is $x\in M$ (in consequence $(x,\dot x)\in TM$).
Observe that the Lagrangians of the type $l: TM\times {\mathfrak
  g}\rightarrow \R$ mainly appears as reduction of Lagrangians of the type $L: T(M\times G)\to \R$, which are invariant under the action of the  Lie group $G$.
Under the identification $T(M\times G)/G\equiv TM\times {\mathfrak g}$
we obtain the reduced Lagrangian $l$. We first develop the discrete optimal
control problem for systems in an unconstrained principle bundle
setting in~\S\ref{sec:docpb}. Nonholonomic constraints are then added
to treat the more general case of locomotion systems
in~\S\ref{sec:docns}.

\subsection{Discrete Optimal Control on Principle Bundles}\label{sec:docpb}
The discrete case is modeled by a Lagrangian $l_d: M\times M\times G\to \R$ which is  an approximation of the action integral in one time step
\[
l_{d}(x_{k},x_{k+1},W_{k})\simeq\int_{hk}^{h(k+1)}l\lp x(t), \dot{x}(t), \xi(t)\rp\,dt,
\]
where $(x_{k},x_{k+1})\in M\times M$ and $W_{k}\in G$. Again, we make an election for the discrete control forces $f^{\pm}_{k}: M\times M\times G \times U\rightarrow T^{*}M\times\mathfrak{g}^{*}$, where $U\subset\R^{m}$:
\begin{eqnarray*}
&&f^{-}_{k}(x_{k},x_{k+1},W_{k},u_{k}^{-})=\lp \bar{f}_{k}^{-}(x_{k},x_{k+1},W_{k},u_{k}^{-}), \hat{f}_{k}^{-}(x_{k},x_{k+1},W_{k},u_{k}^{-})\rp,\\
&&f^{+}_{k}(x_{k},x_{k+1},W_{k},u_{k}^{+})=
\lp \bar{f}_{k}^{+}(x_{k},x_{k+1},W_{k},u_{k}^{+}), \hat{f}_{k}^{+}(x_{k},x_{k+1},W_{k},u_{k}^{+})\rp,
\end{eqnarray*}
here $f_k^-\in T_{x_k}^{*}M\times\mathfrak{g}^{*}$ and $f_k^+\in T_{x_{k+1}}^{*}M\times\mathfrak{g}^{*}$ (more concretely $\bar{f}_{k}^{-}\in T_{x_k}^{*}M$, $\bar{f}_{k}^{+}\in T_{x_{k+1}}^{*}M$, $\hat{f}_{k}^{-}\in\dal$, $\hat{f}_{k}^{+}\in\dal$).


Similarly to the developments in $\S$ \ref{TQ} and $\S$ \ref{AL} we
can formulate the {\bf discrete Lagrange-D'Alembert principle}:
\begin{eqnarray*}
\delta\sum_{k=0}^{N-1}l_{d}(x_{k},x_{k+1},W_{k})&+&\sum_{k=0}^{N-1}\bra f^{-}_{k},(\delta x_{k},\eta_{k})\ket
\\&+&\sum_{k=0}^{N-1}\bra f_{k}^{+},(\delta x_{k+1},\eta_{k+1})\ket=0,
\end{eqnarray*}
which can be rewritten as
\begin{eqnarray*}
\delta\sum_{k=0}^{N-1}l_{d}(x_{k},x_{k+1},W_{k})&+&\sum_{k=0}^{N-1}\bar{f}_{k}^{-}\delta x_{k}+\sum_{k=0}^{N-1}\bar{f}_{k}^{+}\delta x_{k+1}
\\&+&\sum_{k=0}^{N-1}\bra \hat{f}_{k}^{-},\eta_{k}\ket+\sum_{k=0}^{N-1}\bra \hat{f}_{k}^{+},\eta_{k+1}\ket=0,
\end{eqnarray*}
for all variations $\lc\delta x_{k}\rc_{k=0}^{N}$  with $\delta x_k\in T_{x_k} M$ and $\delta x_{0}=\delta x_{N}=0$; also $\lc\delta W_{k}\rc_{k=0}^{N}$ with $\delta W_{k}\in T_{g_k}G$, such that $\delta W_{k}=-\eta_{k}W_{k}+W_{k}\eta_{k+1}$, being $\lc\eta_{k}\rc_{k=0}^{N}$ a sequence of independent elements of $\mathfrak{g}$ such that $\eta_{0}=\eta_{N}=0$.

Applying variations in the last expression and rearranging the sum, we finally obtain the complete set of {\bf forced discrete Euler-Lagrange equations}:
\begin{eqnarray}
&&\hspace{-1cm}D_{1}l_{d}(x_{k},x_{k+1},W_{k})+D_{2}l_{d}(x_{k-1},x_{k},W_{k-1})+\bar{f}_{k}^{-}+\bar{f}_{k-1}^{+}=0,\label{At1}\\
&&\hspace{-1cm} l_{W_{k-1}}^{*}D_{3}l_{d}(x_{k-1},x_{k},W_{k-1})-
r_{W_{k}}^{*}D_{3}l_{d}(x_{k},x_{k+1},W_{k})+\hat{f}^{-}_{k}+\hat{f}^{+}_{k-1}=0,\label{At2}
\end{eqnarray}
with $k=1, \ldots, N-1$.
Since we are dealing with an optimal control problem, we introduce a  discrete cost function $C_d: M\times G\times M\times U\times U\rightarrow \R$. As in previous cases, our objective is to extremize the following sum
\[
\sum_{k=0}^{N-1}C_d(x_k, W_{k}, x_{k+1}, u_k^-, u_k^+),
\]
subjected to equations (\ref{At1}) and (\ref{At2}). Let us initially restrict our attention to the case of fully actuated systems.
\begin{definition}{{\rm {\bf (Fully actuated discrete system)}}}\label{DefAB}
We say that the discrete mechanical control system is fully actuated if the
mappings
\begin{eqnarray*}
&f^-_k\big|_{(x_0, x_1, W_{1})}: U\to T_{x_0}^*M\times {\mathfrak g}^*, \quad  f^-_k\big|_{(x_0, x_1, W_{1})}(u)=f^-_k(x_0, x_1,W_{1}, u),&\\
&f^+_k\big|_{(x_0, x_1, W_{1})}: U\to T_{x_1}^*M\times {\mathfrak g}^*,  \quad  f^+_k\big|_{(x_0, x_1, W_{1})}(u)=f^+_k(x_0, x_1, W_{1}, u)&
\end{eqnarray*}
are both diffeomorphisms.
\end{definition}
According to equations (\ref{At1}) and (\ref{At2}), we can introduce the momenta by means of the following discrete Legendre transforms:
\begin{eqnarray*}
&&p_{k}=-D_{1}l_{d}(x_{k},x_{k+1},W_{k})-\bar{f}_{k}^{-},\\
&&p_{k+1}=D_{2}l_{d}(x_{k},x_{k+1},W_{k})+\bar{f}_{k}^{+},\\\\
&&\mu_{k}=r_{W_{k}}^{*}D_{3}l_{d}(x_{k},x_{k+1},W_{k})-\hat{f}_{k}^{-},\\
&&\mu_{k+1}=l_{W_{k}}^{*}D_{3}l_{d}(x_{k},x_{k+1},W_{k})+\hat{f}_{k}^{+}.
\end{eqnarray*}
In the fully actuated case, is possible to find the value of all control forces in terms of $x_{k},x_{k+1}, W_{k}, p_{k}, p_{k+1},\mu_{k},\mu_{k+1}$, that is:
\begin{eqnarray}
u_k^-&=&u_k^-(x_k, x_{k+1}, W_{k}, p_k, \mu_k)\label{uk-},\\
u_k^+&=&u_k^+(x_k, x_{k+1}, W_{k}, p_{k+1}, \mu_{k+1})\label{uk+}.
\end{eqnarray}
Replacing (\ref{uk-}) and (\ref{uk+}) into $C_d$, we finally obtain the discrete Lagrangian that completely describes our system:
\[
\mathcal{L}_{d}: T^*M\times\dal\times G\times\dal\times T^*M\longrightarrow \R.
\]
The associated {\bf discrete cost functional} is
\begin{equation}\label{ATCF}
\mathcal{J}_{d}=\sum_{k=0}^{N-1}\mathcal{L}_{d}(x_{k},p_{k},\mu_{k},W_{k}, \mu_{k+1}, x _{k+1},p_{k+1}).
\end{equation}
As usual, we take now variations in (\ref{ATCF}) in order to obtain the discrete Euler-Lagrange equations for our optimal control problem (with some abuse of notation we denote $\hat Q_k=(x_{k},p_{k},\mu_{k},W_{k},\mu_{k+1}, x_{k+1},p_{k+1})$ the whole set of coordinates in the new phase space):
\begin{eqnarray*}\label{ELC1-1}
D_6{\mathcal L}_{d}(\hat Q_{k-1})&+&D_1{\mathcal L}_{d}(\hat Q_k)=0\; ,\\\label{ELC2-2}
D_7{\mathcal L}_{d}(\hat Q_{k-1})&+&D_2{\mathcal L}_{d}(\hat Q_k)=0\; ,\\
D_5{\mathcal L}_{d}(\hat Q_{k-1})&+&D_3{\mathcal L}_{d}(\hat Q_k)=0\; ,\label{ELC1-3}\\
l_{W_{k-1}}^* D_4{\mathcal L}_{d}(\hat Q_{k-1})&-&r_{W_{k}}^*D_4{\mathcal L}_{d}(\hat Q_k)=0,\label{ELC1-4}
\end{eqnarray*}
together with the forced discrete Euler-Lagrange equations (\ref{At1})
and (\ref{At2}).

Typically, actuation is achieved by controlling only a subset of the
shape variables. In our setting this is can be regarded as {\bf
  underactuation} -- the mappings in definition \ref{DefAB} become
  embeddings. If this is the case, it is necessary to introduce
constraints and apply constrained variational calculus as in $\S$
\ref{UD1} and $\S$ \ref{AL}.

\subsection{Discrete Optimal Control of Nonholonomic Systems}\label{sec:docns}

This subsection is devoted to add nonholonomic constraints to the picture. Holonomic constraints might be considered as a pacticular case of the nonholonomic ones (see \cite{Sigrid} for further details).  With this extension it would be possible consider examples of optimal control of robotic vehicles. In the following we will expose the theoretical framework, leaving for future research the application to concrete examples.

A  controlled discrete nonholonomic system
 on $M\times M\times G$ is given by the following quadruple (see \cite{IMMM,KFM}):

\begin{itemize}
\item[$i)$]  A {\bf regular discrete Lagrangian} $l_d: M\times M\times G
\rightarrow \R$.

 \item[$ii)$] A {\bf discrete constraint embedded
submanifold} ${\mathcal M}_c$ of $M\times M\times G$.

\item[$iii)$] A {\bf constraint distribution},
${\mathcal D}_c$, which is a vector subbundle  of the vector bundle
$\tau_{_{TM\times\al}}:TM\times {\mathfrak g}\rightarrow M$,  such that $\dim
{\mathcal M}_c = \dim {\mathcal D}_c$. Typically, there is a relation between the constraint distribution and the discrete constraint, since from ${\mathcal M}_c$ we induce
for
every $x\in M$,  the subspace ${\mathcal D}_c(x)$
of $T_x M\times {\mathfrak g}$ given by
\[
{\mathcal D}_c(x)=T_{(x,x,e)}{\mathcal M}_c\cap \lp T_x M\times {\mathfrak g}\rp,
\]
where we are identifying $T_x M\times {\mathfrak g}\equiv 0_x\times
T_{x}M\times T_eG$, with $e$ being the identity element of the Lie group $G$.
\item[$iv)$] The discrete control forces $f^{\pm}_{k}: {\mathcal M}_c \times U\rightarrow T^*M\times {\mathfrak g}^*$ where $U\subset\R^{m}$ (again, forces $f_k^{\pm}$ split into $\bar f_k^{\pm}$ and $\hat f_k^{\pm}$ as in the previous section).
\end{itemize}

We have the following  {\bf discrete version of the Lagrange-D'Alembert principle for controlled nonholonomic systems}:
\begin{eqnarray*}
\delta\sum_{k=0}^{N-1}l_{d}(x_{k},x_{k+1},W_{k})&+&\sum_{k=0}^{N-1}\bra f^{-}_{k},(\delta x_{k},\eta_{k})\ket
\\&+&\sum_{k=0}^{N-1}\bra f_{k}^{+},(\delta x_{k+1},\eta_{k+1})\ket=0,
\end{eqnarray*}
for all variations $\lc\delta x_{k}\rc_{k=0}^{N}$, with $\delta x_{0}=\delta x_{N}=0$; and $\lc\delta W_{k}\rc_{k=0}^{N}$, such that $\delta W_{k}=-\eta_{k}W_{k}+W_{k}\eta_{k+1}$, being $\lc\eta_{k}\rc_{k=0}^{N}$, verifying $(\delta x_k, \eta_k)\in {\mathcal D}_c(x_k)\subseteq T_{x_k}M\times {\mathfrak g}$  such that $\eta_{0}=\eta_{N}=0$.
Moreover, $(x_k, x_{k+1}, W_{k})\in {\mathcal M}_c$, $k=0,\ldots, N-1$ (see \cite{IMMM}).

Take a basis of sections $\{(X^{a}, \tilde{\eta}^{a})\}$ of the vector bundle $\tau_{{\mathcal D}_c}: {\mathcal D}_c\longrightarrow M$, where
$X^a\in\mathfrak{X}(M)$ and $\tilde\eta^a:M\Flder\al$ for
$a=1,...,\text{rank}({\mathcal D}_c)$.  Hence, the equations of motion derived from the discrete Lagrange-D'Alembert principle for controlled nonholonomic systems are:
\begin{small}
\begin{eqnarray}
&& 0=\bra D_{1}l_{d}(x_{k},x_{k+1},W_{k})+ D_{2}l_{d}(x_{k-1},x_{k},W_{k-1}) + \bar{f}_{k}^{-}+ \bar{f}_{k-1}^{+}\,,\, X^{a}(x_k)\ket\nonumber \\
\label{eq-noh}\\&&+ \bra l_{W_{k-1}}^{*}D_{3}l_{d}(x_{k-1},x_{k},W_{k-1})-
r_{W_{k}}^{*}D_{3}l_{d}(x_{k},x_{k+1},W_{k})+\hat{f}^{-}_{k}+\hat{f}^{+}_{k-1}\,,\, \tilde{\eta}^{a}(x_k)\ket,\nonumber\\\nonumber\\
&&0=\Psi^{\alpha}(x_{k},x_{k+1},W_{k}),\label{eq-cons}
\end{eqnarray}
\end{small}
where $\Psi^{\alpha}(x_{k},x_{k+1},W_{k})=0$ are the constraints which locally determine ${\mathcal M}_d$.

In a more geometric way, we can write equations (\ref{eq-noh}) and (\ref{eq-cons}) as follows
\begin{eqnarray*}
0=(i_{{\mathcal D}_c})^* \Big( D_{1}l_{d}(x_{k},x_{k+1},W_{k})&+& D_{2}l_{d}(x_{k-1},x_{k},W_{k-1})+ \bar{f}_{k}^{-}+ \bar{f}_{k-1}^{+}, \\
l_{W_{k-1}}^{*}D_{3}l_{d}(x_{k-1},x_{k},W_{k-1})&-&
r_{W_{k}}^{*}D_{3}l_{d}(x_{k},x_{k+1},W_{k})+\hat{f}^{-}_{k}+\hat{f}^{+}_{k-1}\Big),
\end{eqnarray*}
where$(x_k, x_{k+1}, W_{k})\in {\mathcal M}_c$ and $i_{{\mathcal D}_c}: {\mathcal D}_c\hookrightarrow TM\times {\mathfrak g}$ is the canonical inclusion.

Given a discrete cost function  $C_d: U\times {\mathcal
  M}_c\times U\longrightarrow \R$ and the optimal control problem is
to minimize the action sum
\[
\sum_{k=0}^{N-1} C_d(u_k^-, x_k, W_{k}, x_{k+1}, u_k^+)
\]
subject to equations (\ref{eq-noh}) and (\ref{eq-cons}) and to some
given boundary conditions. We next distinguish between the fully
and under--actuated case using the following definition:
\begin{definition}{{\rm {\bf (Fully actuated nonholonomic discrete system})}}
We say that the discrete nonholonomic mechanical control system is fully actuated if the
mappings
\begin{eqnarray*}
&F^-_k\big|_{(x_0, x_1, W_{1})}: U\to {\mathcal D}_c^*, \quad  F^-_k\big|_{(x_0, x_1, W_{1})}(u)=(i_{{\mathcal D}_c})^*(f^-_k(x_0, x_1, W_{1}, u)),&\\
&F^+_k\big|_{(x_0, x_1, W_{1})}: U\to {\mathcal D}_c^*,  \quad  F^+_k\big|_{(x_0, x_1, W_{1})}(u)=(i_{{\mathcal D}_c})^*(f^+_k(x_0, x_1, W_{1}, u)),&
\end{eqnarray*}
are both diffeomorphisms for all $(x_0, x_1, W_{1})\in {\mathcal M}_c$.
\end{definition}
Regarding equation (\ref{eq-noh}) and its geometric redefinition just below, let introduce the following  momenta:
\begin{eqnarray*}
&&\pi_{k}=(i_{{\mathcal D}_c})^*\left(-D_{1}l_{d}(x_{k},x_{k+1},W_{k})-\bar f_k^-, r_{W_{k}}^{*}D_{3}l_{d}(x_{k},x_{k+1},W_{k})-\hat f _{k}^{-}\right),\\
&&\pi_{k+1}=(i_{{\mathcal D}_c})^*\left(D_{2}l_{d}(x_{k},x_{k+1},W_{k})+\bar f_k^+, l_{W_{k}}^{*}D_{3}l_{d}(x_{k},x_{k+1},W_{k})+\hat {f}_{k}^{+}\right),
\end{eqnarray*}
where both $\pi_k$ and $\pi_{k+1}$ belong to $\mathcal{D}_{c}^{*}$. In the fully actuated case, the value of all control forces can be completely determined in terms of $x_{k},x_{k+1}, W_{k}, \pi_{k}, \pi_{k+1}$, where the coordinates $(x_k, x_{k+1}, W_{k})$ always belong to ${\mathcal M}_c$.
Therefore we can re-express the cost function in terms of these variables and, in consequence,  derive
the discrete Lagrangian
\[
{\mathcal L}_d: \lp{\mathcal D}_c^*\rp
\fourIdx{}{\tau_{{\mathcal D}_c^*}}{}{pr_1}{\times}
 \lp{\mathcal M}_c\rp
\fourIdx{}{pr_2}{}{\tau_{{\mathcal D}_c}^*}{\times}
 \lp{\mathcal D}_c^*\rp
 \rightarrow \R,
\]
where $pr_i: {\mathcal M}_d\subseteq M\times M\times G\rightarrow M$ are the projections onto the first and second arguments and
$\tau_{{\mathcal D}_c^*}: {\mathcal D}_c^*\rightarrow M$ the vector bundle projection.

Observe that we can consider this case as a constrained discrete variational problem taking an extension
\[
\widetilde{{\mathcal L}_d}: {\mathcal D}_c^*\times G \times {\mathcal D}_c^*\rightarrow \R
\]
of ${\mathcal L}_d$ subjected to the constraints $\Psi^{\alpha}(x_{k},x_{k+1},W_{k})=0$.

Therefore, denoting $\hat Q_k=(x_{k},\pi_{k},W_{k},x_{k+1},\pi_{k+1})$ as the whole set of coordinates of the new phase space ${\mathcal D}_c^*\times G \times {\mathcal D}_c^*$, we deduce that the equations of motion are
\begin{eqnarray*}
D_4\widetilde{{\mathcal L}_{d}}(\hat Q_{k-1})+D_1\widetilde{{\mathcal L}_{d}}(\hat Q_k)&=&\lambda^{k-1}_{\alpha} D_2\Psi^{\alpha}(x_{k-1}, x_{k}, W_{k-1})\nonumber
\\
&&+\lambda^k_{\alpha} D_1\Psi^{\alpha}(x_{k}, x_{k+1}, W_{k}),\; \\\\\label{ELC2-noh-2}
D_5\widetilde{{\mathcal L}_{d}}(\hat Q_{k-1})+D_2\widetilde{{\mathcal L}_{d}}(\hat Q_k)&=&0\;  ,\label{ELC1-noh-3}\\\\
l_{W_{k-1}}^* D_3\widetilde{{\mathcal L}_{d}}(\hat Q_{k-1})-r_{W_{k}}^*D_3\widetilde{{\mathcal L}_{d}}(\hat Q_k)&=&
\lambda^{k-1}_{\alpha}l_{W_{k-1}}^* D_3\Psi^{\alpha}(x_{k-1}, x_{k}, W_{k-1})\nonumber\\
&&-\lambda^{k}_{\alpha}r_{W_{k}}^*D_3\Psi^{\alpha}(x_{k}, x_{k+1}, W_{k}),
\label{ELC1-noh4}\\\\
\Psi^{\alpha}(x_{k},x_{k+1},W_{k})&=&0\; , \label{ELC1-noh5}
\end{eqnarray*}
where $\lambda_{\alpha}^k$ are the Lagrange multipliers of the new
constrained problem. The underactuated case can be handled by adding
new constraints and applying discrete constrained variational calculus
similarly to~\S\ref{RED}.

A natural framework that simplifies the previous construction is based on discrete mechanics on Lie groupoids \cite{groupoid}. The Lie groupoid structure generalizes the case of $Q\times Q$, the Lie group $G$ and also many intermediate situations. In particular, many of the examples studied in this paper can be modeled using Lie groupoid techniques adapted to our formalism (see \cite{JMdD}).

\section{Conclusions}
This paper develops numerical methods for optimal control of
Lagrangian mechanical systems defined on tangent bundles, Lie groups,
trivial principal bundles, and nonholonomic systems. The proposed
approach preserves the geometry and
variational structure of mechanics through the discretization of
the variational principles on manifolds. The key point is to solve
the optimal control through discrete mechanics, i.e. by formulating the
optimization as the solution of an action principle of
a higher-dimensional system in a new Lagrangian phase space,
i.e. $T^*Q\times T^*Q$ in the general case and $\dal\times
G\times\dal$ in the Lie group case. The optimal control algorithm is
then derived as a variational integrator subject to boundary
conditions. We thus expect that both the dynamics and optimal control
solutions will have accurate and stable numerical behavior (due to
symplectic-momentum preservation) even at large
time-steps (which allows for improved run-time efficiency).



Simulations of an underactuated underwater vehicle illustrate an
application of the method. Yet, further numerical studies and
comparisons would be necessary to exactly quantify the advantages and
the limitations of the proposed algorithm.
An important future direction is thus to study the convergence
properties of the optimal control system. Convergence for general
nonlinear systems is a complex issue. In this respect, it is
interesting to note that the discrete mechanics and optimal control on
Lie groups such as the example in~\label{Vehiculo} using the Cayley
map results in polynomial form without further approximation or Taylor
series truncation. A useful future direction is then to study the
regions of attraction of the numerical continuation using tools from
algebraic geometry. 

More generally, the theoretical framework introduced
in~\S\ref{sec:ext} can serve as a basis for deriving algorithms for
control systems such as multi-body locomotion systems or robotic
vehicles with nonholonomic
constraints. Furthermore, the developed classes of systems can be
unified through the recently developed groupoid
framework~\cite{IMMM,weinstein}. Each of the considered product spaces
(e.g. $Q\times Q$) can be regarded as a single groupoid space with
equations of motion resulting from a single generalized discrete
variational principle. This will enable the automatic solution of
optimal control problems for various complex systems and a convenient
unified framework for implementing practical optimization schemes such
as~\cite{BlHuLeSa2009,KM,LeMcLe2006,objuma}. More importantly, this
viewpoint can be used to apply standard discrete Lagrangian regularity
conditions (e.g.~\cite{mawest}) to optimal control problems evolving on
the groupoid space. This would provide a deeper insight into the
solvability of the resulting optimization schemes.


\section*{Appendix A: Lemmae}

\begin{lemma}\label{Ant}
(see \cite{MaRa})
Let $g\in G$, $\lambda\in\al$ and $\delta f$ denote the variation of a function $f$ with respect to its parameters. Assuming $\lambda$ is constant, the following identity holds
{\rm\[
\delta(\Ad_{g}\,\lambda)=-\Ad_{g}\,[\lambda\,,\,g^{-1}\delta g],
\]}
where $[\cdot\,,\,\cdot]:\al\times\al\Flder\R$ denotes the Lie bracket operating or equivalently $[\xi\,,\,\eta]\equiv\ad_{\xi}\eta$, for given $\eta,\,\xi\in\al$.
\end{lemma}



\begin{lemma}(see \cite{Rabee})\label{lemaPr}
  The following identity holds
 {\rm \[
  \mbox{d}\tau_{\xi}\,\eta=\Ad_{\tau(\xi)}\,\mbox{d}\tau_{-\xi}\,\eta,
  \]}
  for any $\xi,\eta\in\al$.
\end{lemma}

\begin{lemma}(see \cite{Rabee})\label{lemaPr2}
  The following identity holds
 {\rm \[
  \mbox{d}\tau^{-1}_{\xi}\,\eta=\mbox{d}\tau_{-\xi}^{-1}\lp\Ad_{\tau(-\xi)}\,\eta\rp,
  \]}
  for any $\xi, \eta \in\al$.
\end{lemma}

\end{document}